\documentclass[onefignum,onetabnum]{siamonline190516}

\usepackage{calrsfs,amsmath,amssymb,graphicx,array,accents,figlatex,epsfig}
\usepackage{setspace,wrapfig,colortbl,marvosym, bbm, stmaryrd,centernot, color}
\usepackage{amsfonts}
\usepackage{epstopdf}
\usepackage{algorithmic}
\ifpdf
  \DeclareGraphicsExtensions{.eps,.pdf,.png,.jpg}
\else
  \DeclareGraphicsExtensions{.eps}
\fi

\usepackage{enumitem}
\setlist[enumerate]{leftmargin=.5in}
\setlist[itemize]{leftmargin=.5in}

\newsiamremark{rem}{Remark}
\newsiamthm{assumption}{Assumption}
\newsiamremark{remark}{Remark}
\newsiamremark{hypothesis}{Hypothesis}
\crefname{hypothesis}{Hypothesis}{Hypotheses}
\newsiamthm{claim}{Claim}

\newcommand{\R}{\mathbb{R}}

\headers{Differential tomography of micromechanical evolution}{F. Pourahmadian and H. Haddar}

\title{Differential tomography of micromechanical evolution in elastic materials of unknown micro/macrostructure\thanks{
\funding{This work was funded by the University of Colorado Boulder through FP's startup and IS-IRT grant.}}}

\author{Fatemeh Pourahmadian\thanks{Department of Civil, Environmental and Architectural Engineering and Department of Applied Mathematics, University of Colorado Boulder, USA 
  (\email{fatemeh.pourahmadian@colorado.edu)}.}
\and Houssem Haddar\thanks{INRIA Saclay Ile de France and Ecole Polytechnique (CMAP) Route de Saclay, F-91128, Palaiseau, France 
  (\email{houssem.haddar@polytechnique.edu)}.}}

\usepackage{amsopn}

\ifpdf
\hypersetup{
  pdftitle={Differential tomography of micromechanical evolution in elastic domains of unknown micro/macrostructure},
  pdfauthor={F. Pourahmadian and H. Haddar}
}
\fi

\newlength{\kaka}

\newcommand{\ahref}[2]{}

\newcommand{\Zsup}{^{\text{\tiny o}}}

\newcommand{\obs}{^{\text{obs}}}

\newcommand{\ff}{^{\text{\tiny f}}}

\newcommand{\beq}{\begin{equation}}
\newcommand{\eeq}{\end{equation}}
\newcommand{\lb}{\label}

\newcommand{\bea}{\begin{eqnarray}}
\newcommand{\eea}{\end{eqnarray}}
\newcommand{\bxr}{\begin{array}}
\newcommand{\exr}{\end{array}}

\newcommand\exs{\hspace*{0.4mm}}

\newcommand\nxs{\hspace*{-0.2mm}}

\newcommand{\norms}[1]{\parallel\! #1 \!\parallel}

\newcommand{\bSig} {\boldsymbol{\Sigma}}
\newcommand{\bW} {\boldsymbol{W}}
\newcommand{\bC} {\boldsymbol{C}}
\newcommand{\bK} {\boldsymbol{K}}

\newcommand{\bR} {\boldsymbol{\sf R}}

\newcommand{\bn} {\boldsymbol{n}}
\newcommand{\ba} {\boldsymbol{a}}
\newcommand{\bb} {\boldsymbol{b}}

\newcommand{\bx} {\boldsymbol{x}}
\newcommand{\by} {\boldsymbol{y}}
\newcommand{\be} {\boldsymbol{e}}

\newcommand{\bg} {{\boldsymbol{g}}}
\newcommand{\bq} {{\boldsymbol{q}}}

\newcommand{\bfT} {\boldsymbol{T}}

\newcommand{\bd} {\boldsymbol{d}}
\newcommand{\bI} {\boldsymbol{I}}

\newcommand{\pff}{\btu}

\newcommand{\sip} {\!\cdot\!}

\newcommand{\OO}{\Omega}

\newcommand{\OOd}{{\Omega}}

\newcommand{\bzero}{\boldsymbol{0}}

\newcommand{\bu} {\boldsymbol{u}}

\newcommand{\bt} {{\boldsymbol{t}}}

\newcommand{\bv} {\boldsymbol{v}}

\newcommand{\bxi} {\boldsymbol{\xi}}
\newcommand{\bxio} {\bxi\Zsup}

\newcommand{\dbvz} {\llbracket {\bv_{\!\circ}}\rrbracket}
\newcommand{\dbvk} {\llbracket {\bv_\kappa}\rrbracket}
\newcommand{\dbva} {\llbracket {\bv_\alpha}\rrbracket}
\newcommand{\dbvap} {\llbracket {\bv_{\alpha+1}}\rrbracket}

\newcommand{\btu} {\text{\bf{u}}}

\newcommand{\bPhi}{\boldsymbol{\Phi}}

\newcommand{\bphi}{{\boldsymbol{\varphi}}}
\newcommand{\bpsi}{{\boldsymbol{\psi}}}

\begin{document}

\maketitle

\begin{abstract}
\emph{Differential evolution indicators} are introduced for 3D spatiotemporal imaging of micromechanical processes in complex materials where progressive variations due to manufacturing and/or aging are housed in a highly scattering background of a-priori unknown or uncertain structure. In this vein, a three-tier imaging platform is established where:~(1)~the domain is periodically (or continuously) subject to illumination and sensing in an arbitrary configuration;~(2)~sequential sets of measured data are deployed to distill segment-wise scattering signatures of the domain's internal structure through carefully constructed, non-iterative solutions to the scattering equation; and~(3)~the resulting solution sequence is then used to rigorously construct an imaging functional carrying appropriate invariance with respect to the unknown stationary components of the background e.g.,~pre-existing interstitial boundaries and bubbles. This gives birth to differential indicators that specifically recover the 3D support of micromechanical evolution within a network of unknown scatterers. The direct scattering problem is formulated in the frequency domain where the background is comprised of a random distribution of monolithic fragments. The constituents are connected via highly heterogeneous interfaces of unknown elasticity and dissipation which are subject to spatiotemporal evolution. The support of internal boundaries are sequentially illuminated by a set of incident waves and thus-induced scattered fields are captured over a generic observation surface. The performance of the proposed imaging indicator is illustrated through a set of numerical experiments for spatiotemporal reconstruction of progressive damage zones featuring randomly distributed cracks and bubbles.
\end{abstract}



\section{Introduction} \label{sec1}

Fast waveform tomography solutions germane to uncertain (\nxs or unknown\nxs) environments bear direct relevance to (a) timely detection of degradation in safety-sensitive components~\cite{cobl2015,barn2013}, and (b) in-situ monitoring of additive manufacturing (AM) processes~\cite{ever2016}. In nuclear power plants, for instance, critical components such as reactor and fuel cells are comprised of composite materials whose topology and properties are uncertain at micro-, meso-, and macro-scales as a result of manufacturing and/or aging. The deterioration of these materials due to various chemo-physical mechanisms such as irradiation and thermal cycling are not yet fully understood~\cite{Naus1996}. These processes, however, spur continuous microstructural evolution leading to an inevitable development of anomalies responsible for the loss of structural integrity and diminished functional performance. Thus, timely detection of deterioration at the microstructure scale and active spatiotemporal tracking of their evolution are paramount for early and robust mitigation of damage in such systems. In advanced manufacturing, one of the main challenges is online evaluation of the AM performance~\cite{fraz2014,hua2015}, demanding real-time in-situ characterization of components during fabrication~\cite{tahe2018}. In this vein, a major hindrance is that the scattering signature of evolving regions is often eluded by the footprints of unknown stationary scatterers in a complex specimen. A sensing scheme amenable to highly scattering environments will contribute to (a) better understanding of the manufacturing process and its implications on the quality of the final product, and (ii) optimal design and closed-loop control of the AM processes. 

Recently, major progress has been made toward developing robust imaging solutions that enable real-time sensing in complex materials~\cite{rose2014,matl2015,roja2015,amin2017,cobl2015,hung2013,hank2016,Guzi2014}. State-of-the-art examples include:~ultrasonic surface wave methods~\cite{rose2014}, nonlinear ultrasound~\cite{matl2015}, penetrating-radar techniques~\cite{amin2017}, infrared thermography~\cite{cobl2015}, laser shearography~\cite{hung2013}, X-ray computed tomography~\cite{hank2016} and acoustic tomography imaging~\cite{Guzi2014}. So far, these developments mostly rely upon (a) simplistic characterization of the background disregarding uncertain yet fundamental features (such as interstitial boundaries) across multiple scales~\cite{tokm2013},~{(b)}~significant assumptions on the nature of wave motion~\cite{rose2014},~{(c)}~partial data inversion deploying only a few signatures of the scattered field measurements~\cite{roja2015}, and~{(d)}~data processing schemes amenable to ad hoc sensing configurations~\cite{cobl2015}. Such attributes impose a number of limitations, including:~{(i)}~high sensitivity to the assumed background structure,~{(ii)}~insensitivity to less-understood properties of microstructured materials,~{(iii)}~major restrictions in terms of the location of ultrasonic transducers,~{(iv)}~limited scalability beyond laboratory applications e.g., for the purpose of in-situ monitoring. Therefore, there is a critical need for the next generation of imaging tools that transcend some of these barriers. 

Ongoing efforts in this vein are mainly focused on optimization-based approaches to waveform inversion. These technologies typically incur high computational cost as a crucial obstacle to real-time sensing, and high sensitivity to unknown features of the background leading to multiple sets of ``optimal'' solutions and thus, ambiguity of the results. More recently, approaches to fast waveform tomography~\cite{Audibert2014,audi2017,bonn2019,audi2015,cako2016,de2018} have been brought under the spotlight for their capabilities pertinent to imaging in highly scattering media. While this class of inverse solutions generally demand an a-priori characterization of the background for their successful performance, most recent developments including the present study indicate that this requirement could be relaxed in presence of sequential measurements, generating a suit of new imaging modalities that surpass some of the existing limitations.

In particular, this study takes advantage of some recent advances in design of sampling methods~\cite{Fatemeh2017(2),audi2015, had2017, Fatemeh2017} to develop a non-iterative full-waveform approach for spatiotemporal tracking of progressive variations in complex materials. The idea is to deploy sequential sets of scattered field measurements to rigorously construct an imaging functional endowed with appropriate invariance with respect to (unknown) stationary components of the background e.g., its time-invariant scatterers. The resulting differential indicators uniquely characterize the support of evolution without the need to reconstruct the entire domain across pertinent scales which may be practically insurmountable. In the case of volumetric scatterers, such invariants are furnished via solutions to the so-called interior transmission problem~and the relation between two such solutions before and after the evolution~\cite{audiDLSM,cako2016}. The key observation in developing the imaging functional is that such solutions may be approximated by using sampling type techniques~\cite{audi2015, Fatemeh2017}. 

The fundamental challenge which impedes direct extension of these advances to ultrasonic imaging is the existence of non-volumetric scatterers in solid-state materials e.g., interstitial boundaries, fractures, and dislocation networks~\cite{de2018,pour2018(2)}. So that functionals of desired invariance may not be established through the analysis of elastodynamic interior transmission problems. This work aims to address this challenge by studying imaging functionals pertinent to elastic backgrounds with random interfaces and discontinuities across scales. Our analysis is based on the boundary integral representation of scattering solutions, enabling rigorous formulation of invariant quantities critical for establishing differential imaging functionals for such media. The designed indicators are then synthetically tested and validated in a few example configurations featuring randomly distributed interfaces and bubbles. 

This paper is organized as follows. \cref{PS} presents the direct scattering formulation and admissibility conditions on interstitial boundaries so that the forward problem is wellposed. \cref{Prelim} defines the scattering operator and briefly recaps some known results on the properties of this operator for later reference. The differential evolution indicators are introduced and analyzed in \cref{DEI}. \cref{IR} is dedicated to numerical implementation and validation of this imaging solution.

\section{Problem statement}\label{PS}

\begin{figure}[h!]
\begin{center} 
\includegraphics[width=0.99\linewidth]{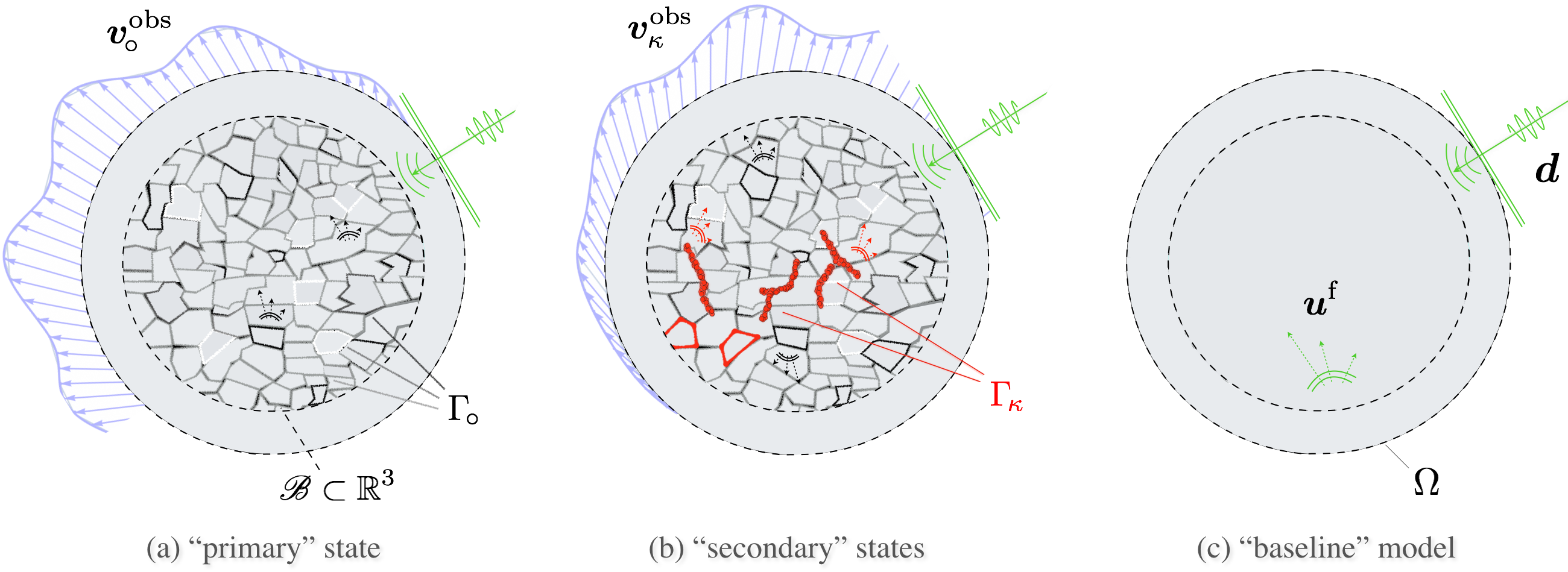}
\end{center} \vspace*{-2.5mm}
\caption{\small{Elastic-wave sensing of microstructural evolution in a background $\mathcal{B} \subset \R^3$ featuring a random network of pre-exiting interfaces $\Gamma_{\!\circ}$: {\bf(a)} \underline{\emph{primary experiments}} conducted at $t \!=\! t_{\nxs\circ}$ via a set of P- and S- plane waves propagating in direction $\bd \in \Omega$, inducing the incident field $\bu^{\textsf{f}}$ in the baseline system shown in (c); the action of $\bu^{\textsf{f}}$ on $\Gamma_{\!\circ}$ results in the scattered field $\bv_{\!\circ}^{\textrm{obs}}$ captured at the far field over the unit sphere of observation angles $\Omega$,~{\bf(b)}~\underline{\emph{secondary experiments}} performed in a similar setting at a sequence of time steps $t_\kappa \!=\! \lbrace t_1, t_2, ... \rbrace$ -- when active chemo-physical processes have created new (or evolved) interfaces $\Gamma_{\!\kappa}$ in $\mathcal{B}$, resulting in the scattered field measurements $\bv_{\kappa}^{\textrm{obs}}$, and~{\bf(c)}~\underline{\emph{baseline model}} of the system constructed synthetically based on a-priori available knowledge of the domain.}}
\label{comp_ps}\vspace*{-4.0mm}
\end{figure}

With reference to~\cref{comp_ps}, consider sequential illumination of microstructural evolution in an elastic domain at sensing steps $t = \lbrace t_{\circ}, t_1, t_2, ... \rbrace$. The domain $\mathcal{B} \!\subset\! \R^3$ at the outset of sensing $t \!=\! t_\circ\!$ is comprised of a random distribution of monolithic fragments of Lipschitz support with mass density $\rho$ and Lam\'{e} parameters $\mu$ and~$\lambda$, connected via perfect or imperfect interfaces $\Gamma_{\!\circ} \!\subset\! \mathcal{B}$. The contact condition at the surface of $\Gamma_{\!\circ}$ is characterized by a symmetric, complex-valued and heterogeneous interfacial stiffness matrix $\bK_{\!\circ}(\bxi), \, \bxi \in \Gamma_{\!\circ}$. Internal interfaces are subject to variations e.g., driven by chemo-physical reactions so that at any secondary sensing step $t = t_\kappa \!\nxs>\! t_\circ$, the domain features newborn and evolved interfaces  $\Gamma_\kappa \subset \mathcal{B}$ endowed with the contact stiffness $\bK_{\!\kappa}(\bxi), \, \bxi \in \Gamma_{\!\kappa}$.  

\begin{assumption}\label{tfc}
Let us denote by $\Gamma_{\!\! N}^t$ the support of all traction-free cracks in $\mathcal{B}$ at time $t$ i.e., $\bK_{\! t} = \bzero$ on $\Gamma_{\!\! N}^t$. In this study, it is assumed that (a) no subset of $\Gamma_N^t$ constitutes a closed surface, and (b) $\mathcal{B} \backslash \Gamma_{\!\! N}^t$ remains  connected.
\end{assumption}

 Let $\Omega$ denote the unit sphere centered at the origin. Given a propagation direction $\bd\in\Omega$ and polarization amplitudes $\bq_p,\bq_s\!\in\mathbb{R}^3$ where $\bq_p\!\parallel\bd$ and~$\bq_s\!\perp\!\bd$, the domain $\mathcal{B}(t)$ is illuminated at every sensing step via a combination of plane P-~and S-~waves so that the incident field takes the form 
\vspace*{-1mm}
\beq\lb{plwa}
\bu\ff(\bxi) ~=~ \bq_p \exs e^{\textrm{i} k_p \bxi \cdot \bd} \:+\: \bq_s \exs e^{\textrm{i} k_s \bxi \cdot \bd}, \qquad \bd \in \Omega, \,\, \bxi \in \R^3,
\eeq
where $k_p$ and $k_s=k_p\sqrt{(\lambda\!+\!2\mu)/\mu}$ denote the respective wave numbers affiliated with the ``baseline" system shown in~\cref{comp_ps}~(c).  At the primary sensing step $t= t_\circ$, the interaction of $\bu\ff$ with the scatterers $\Gamma_{\!\circ}(\bxi)$ gives rise to the scattered field $\bv_{\!\circ}\in H^1_{\mathrm{loc}}(\R^3\backslash\Gamma_{\!\circ})^3$ solving 
\vspace*{-1mm}
\beq\lb{GE}
\begin{aligned}
&\nabla \sip (\bC \colon \! \nabla \bv_{\!\circ}) \,+\, \rho \exs \omega^2\bv_{\!\circ} ~=~ \bzero \quad &\text{in}& \quad {\R^3}\backslash\Gamma_{\!\circ}, \\*[0.5mm]
&\bn_{\nxs\circ} \!\cdot\nxs \bC \exs \colon \!  \nabla  \bv_{\!\circ}~=~ \bK_{\!\circ}(\bxi) \dbvz  \,-\, \bt\ff_{\circ}  \quad &\text{on}& \quad \Gamma_{\!\circ},
\end{aligned}     
\eeq
where $\omega^2=k_s^2 \mu/\rho$ is the frequency of excitation; $\dbvz=[\bv_{\nxs\circ}^+\!-\bv_{\nxs\circ}^-]$ is the jump in~$\bv_{\!\circ}$ across~$\Gamma_{\!\circ}$, hereon referred to as the fracture opening displacement \textcolor{black}{(FOD)}; 
\vspace*{-2mm}
\beq\label{bC}
\bC \:=\: \lambda\,\bI_2\!\otimes\bI_2 \:+\: 2\mu\,\bI_4 \vspace*{-2mm}
\eeq
signifies the fourth-order elasticity tensor; $\bI_m \,(m\!=\!2,4)$ represents the $m$th-order symmetric identity tensor; \mbox{$\bt\ff_{\nxs\circ} = \bn_{\nxs\circ} \cdot \bC \colon \! \nabla \bu\ff$} is the free-field traction vector; $\bn_{\nxs\circ} = \bn_{\nxs\circ}^-$ is the unit normal on~$\Gamma_{\!\circ}$. At subsequent sensing steps $t \!=\! t_\kappa$, the interaction of $\bu\ff$ with $\Gamma_{\!\circ} \cup \Gamma_\kappa$ results in the scattered field $\bv_\kappa \in H^1_{\mathrm{loc}}({\R^3}\backslash\lbrace\Gamma_{\!\circ} \cup \Gamma_\kappa \rbrace)^3$ satisfying
\vspace*{-1.5mm}
\beq\lb{GEi}
\begin{aligned}
&\nabla \sip (\bC \colon \! \nabla \bv_\kappa) \,+\, \rho \exs \omega^2\bv_\kappa ~=~ \bzero \quad &\text{in}& \quad {\R^3}\backslash\lbrace\Gamma_{\!\circ} \cup \Gamma_\kappa \rbrace, \\*[0.5mm]
&\bn_\kappa \!\cdot\nxs \bC \exs \colon \!  \nabla  \bv_\kappa~=~ \bK_\kappa(\bxi) \dbvk  \,-\, \bt\ff_\kappa  \quad &\text{on}& \quad \Gamma_\kappa, \\*[0.75mm]
&\bn_{\nxs\circ} \!\cdot\nxs \bC \exs \colon \!  \nabla  \bv_\kappa~=~ \bK_{\!\circ}(\bxi) \dbvk  \,-\, \bt\ff_{\circ}  \quad &\text{on}& \quad \Gamma_{\!\circ},
\end{aligned}     
\vspace*{-1mm}
\eeq
where $\dbvk=[\bv_\kappa^+\!-\bv_\kappa^-]$ denotes FOD across~$\Gamma_\kappa$; \mbox{$\bt\ff_\kappa = \bn_\kappa \cdot \bC \colon \! \nabla \bu\ff$} is the free-field traction over $\Gamma_\kappa(\bxi)$, and $\bn_\kappa = \bn_\kappa^-$ is the unit normal on~$\Gamma_\kappa$. Formulations of the direct scattering problems~\eqref{GE} and~\eqref{GEi} are complete by imposing the Kupradze radiation condition~\cite{Kuprad1979} on $\bv_{\!\circ}$ and $\bv_\kappa$ at far field. More specifically, on uniquely decomposing the scattered fields into irrotational and solenoidal parts as $\bv_\alpha = \bv_{p_\alpha}\nxs \oplus \exs\bv_{s_\alpha}$ for $\alpha \in \lbrace\circ, \kappa\rbrace$ \textcolor{black}{where 
\vspace*{-2mm}
\beq\label{vpvs}
\bv_{p_\alpha} = \frac{1}{k_s^2\!-\!k_p^2}(\Delta+k_s^2)\bv_\alpha, \qquad \bv_{s_\alpha} = \frac{1}{k_p^2\!-\!k_s^2}(\Delta+k_p^2)\bv_\alpha, \qquad \alpha = \lbrace\circ, \kappa\rbrace, \vspace*{-2mm}
\eeq}
the Kupradze condition can be stated as 
\vspace*{-0.5mm}
\beq\lb{KS}
\frac{\partial\bv_{p_\alpha}}{\partial r} - \text{i} k_p \bv_{p_\alpha} = o\big(r^{-1}\big) \quad \mbox{ and } \quad 
\frac{\partial\bv_{s_\alpha}}{\partial r} - \text{i} k_s \bv_{s_\alpha} = o\big(r^{-1}\big) \qquad \text{as} ~~r:=|\bxi|\to\infty,
\vspace*{-1.5mm}
\eeq    
uniformly with respect to $\hat\bxi:=\bxi/r$. The following remarks shine more light on some specific aspects of the ensuing developments.

 {\bfseries\slshape Background domain.}~Here, our primary knowledge of the system is assumed to be at the ``baseline" level, shown in~\cref{comp_ps}~(c), which is simplistic and mostly limited to idealistic design parameters. However, chemo-physical processes of interest such as early-stage degradation mostly reside at the micro- and meso- scales, developing in a network of pre-existing scatterers of similar scale yet uncertain nature. Accordingly, as illustrated in~\cref{comp_ps}~(a), the background is modeled by an elastic domain endowed with arbitrary interstitial boundaries of heterogeneous contact condition -- spanning from perfectly bonded to traction-free interfaces. This provides a versatile platform for a range of micromechanical phenomena e.g., degradation as a generic cloud of (stationary and evolving) micro-scatterers of random distribution.

 {\bfseries\slshape Anatomy of evolution.}~Stress concentration, chemical reaction, and early-stage irradiation are common producers of \emph{interfacial} damage at micro- and meso- scales~\cite{ross2017,maev2016}. Thermal cycling, fatigue, and shock-wave loading, however, are mostly responsible for distributed fracture zones~\cite{live2014}. Thus, active processes in this work are identified, according to~\cref{comp_ps}~(b), with connected or unconnected sets of heterogeneous \emph{fractures}~$\Gamma_{\!\kappa}$ of interfacial elasticity~$\bK_{\nxs\kappa}(\bxi)$.

  {\bfseries\slshape Illumination frequency.}~The proposed differential imaging scheme is rooted in the sampling methods~\cite{kirs008,Bouk2013,Bellis2013} recognized for providing good quality reconstruction of hidden scatterers at resolution scales transcending the traditional limits of NDE. Here, the illuminating wavelength $\lambda_s = 2 \pi / k_s$ is considered to be comparable with the characteristic length scale of the sought-for processes e.g.,~micro-meso-scale features are probed by micro-meso-scale waves. It is worth mentioning that for multiscale characterization, multi-frequency illumination i.e.,~input signals with appropriate spectral content may be adopted~\cite{guzi2010}.

 {\bfseries\slshape Dimensional platform.}~In what follows, all quantities are rendered \emph{dimensionless} by taking $\rho$, $\mu$, and ${\ell_\circ}$ -- denoting the minimum length scale affiliated with internal boundaries, as the respective reference scales for mass density, elastic modulus, and length -- which amounts to setting $\rho = \mu = \ell_\circ = 1$~\cite{Scaling2003}.

\subsection*{Wellposedness of the sequential direct scattering problems}\lb{WelD} 

\begin{assumption}\label{welp}
In this study, it is assumed that $\bt\ff_{\!\circ} \in H^{-1/2}(\Gamma_{\!\circ})^3$ (resp.~$\bt\ff_\kappa \!\in\! H^{-1/2}(\Gamma_\kappa)^3$) and that $\bK_{\!\circ} \!\in\! L^\infty(\Gamma_{\!\circ})^{3\times 3}$ (resp.~$\bK_{\nxs\kappa} \in L^\infty(\Gamma_\kappa)^{3\times 3}$) is symmetric while satisfying $\overline{\boldsymbol{\theta}} \sip \Im \bK_{\!\circ}(\bxi)\sip \boldsymbol{\theta} \leqslant 0$, $\forall\boldsymbol{\theta}\in\mathbb{C}^3, \, \bxi \in \Gamma_{\!\circ}$ (resp.~$\overline{\boldsymbol{\theta}} \sip \Im \bK_{\nxs\kappa}(\bxi)\sip \boldsymbol{\theta} \leqslant 0$, $\forall\boldsymbol{\theta}\in\mathbb{C}^3, \, \bxi \in \Gamma_\kappa$). 
\end{assumption}

 Under~\cref{tfc} and~\cref{welp}, the direct scattering problems~\cref{GE} and~\cref{GEi} are wellposed. The proof draws from~(a)~Lemma 3.1 and Theorem 3.2 of~\cite{Fatemeh2017} and arguments of unique continuation principles. The proof of~\cite[Theorem 3.2]{Fatemeh2017} directly substantiates that~\cref{GE} and~\cref{GEi} are of Fredholm type, and thus, their wellposedness is certified as soon as the uniqueness of a solution is guaranteed. To verify the latter, let $\bt\ff_{\circ} = \bt\ff_\kappa = \bzero$, then according to~\cite[Lemma 3.1]{Fatemeh2017}, the scattered waveforms $\bv_\alpha(\bxi)$, $\alpha \in \lbrace\circ, \kappa\rbrace$ vanish at the far field as $|\bxi|\to\infty$. The argument is then followed for the case of fragmented backgrounds shown in~\cref{comp_ps}~(a)~and~(b) where~$\mathcal{B}$ is described as a union of simply connected bounded domains $\mathcal{D}_i$, $\lbrace i = 1,2,...,N_{\mathcal{D}} \rbrace$ of Lipschitz boundaries denoted by $\Gamma_{\!\circ}$ (\emph{resp.}~$\Gamma_{\!\circ} \cup \Gamma_\kappa$) at the primary (\emph{resp.}~secondary) sensing step. On setting $\mathcal{D}_0 = \R^3 \!\setminus\! \mathcal{B}$ and $\mathcal{G}_0  = \lbrace 0 \rbrace$, lets define 
\vspace*{-1.5mm}
\beq\lb{ind}
\mathcal{G}_j = \Big{\lbrace}  i \,\,\, \Big{|} \,\,\, i \notin {\textstyle \bigcup\limits_{\kappa =0}^{j-1}} \mathcal{G}_{\kappa} \,\,\, \textrm{and} \,  {\textstyle \bigcup\limits_{\kappa \in \mathcal{G}_{j-1}}} \!\!\!\!\! \mathcal{D}_\kappa  \cap \mathcal{D}_i = {\text{\sf S}}  \Big{\rbrace}, \qquad  i = 1,2,...,N_{\mathcal{D}} , \,\,\,  j = 1,2,...,N_{\text{\tiny $\mathcal{G}$}}, \vspace{-1.2 mm}
\eeq          
where $\sf{S}$ identifies any piecewise analytic surface in $\mathcal{B}$. In light of~\eqref{ind}, the domain may be partitioned into $N_{\text{\tiny $\mathcal{G}$}}+1$ layers $L_j =  \mathop{\cup}_{\kappa \in \mathcal{G}_{j}} \mathcal{D}_\kappa, \, j = 0,1, ..., N_{\text{\tiny $\mathcal{G}$}}$, such that successive application of the unique continuation theorem and Holmgren's principle in each layer completes the uniqueness proof. More specifically, starting from $L_0 = \R^3 \!\setminus\! \mathcal{B}$ where $\bv_\alpha = \bzero$ as $|\bxi|\to\infty$, the unique continuation theorem is deployed to infer $\bv_\alpha = \bzero$ in $L_0$. Subsequently, the jump in displacement $\llbracket \bv_\alpha \rrbracket(\bxi)$ vanishes over the interface of $\bxi \in L_0 \cap L_1 \subset \Gamma_{\!\circ} \cup \Gamma_\kappa$ according to the elastic boundary conditions over $\Gamma_{\!\circ}$ and $\Gamma_\kappa$ in~\cref{GE} and~\cref{GEi}. In this setting, Holmgren's theorem implies that the scattered field $\bv_\alpha$ vanishes in an open neighborhood of $L_0 \cap L_1$ which by virtue of the unique continuation theorem leads to $\bv_\alpha(\bxi) = \bzero$ in $\bxi = L_1$. On repeating this argument in $L_2, ..., L_{\text{\tiny $\mathcal{G}$}}$, one arrives at $\bv_\alpha(\bxi) = \bzero$ in $\mathcal{B} \backslash \lbrace \Gamma_{\!\circ} \cup \Gamma_\kappa \rbrace$ which completes the proof for the uniqueness of a scattering solution in $\mathcal{B}$, and thus, substantiates the wellposedness of the forward problem.      

The scattered waveforms $\bv_\circ \in H^1_{\mathrm{loc}}({\R^3}\backslash\Gamma_{\!\circ})^3\!$ and $\bv_\kappa \in H^1_{\mathrm{loc}}({\R^3}\backslash\lbrace\Gamma_{\!\circ} \cup \Gamma_\kappa \rbrace)^3$ are sequentially captured at $t = \lbrace t_{\circ}, t_1, t_2, ... \rbrace$ in the form of far-field patterns $\bv^{\infty}_{\alpha}\nxs = \bv^{\infty}_{p_\alpha} \nxs\oplus \bv^{\infty}_{s_\alpha}$, $\alpha \in \lbrace\circ, \kappa\rbrace\nxs$, according to the asymptotic expansion
\vspace*{-1mm}
\beq\lb{vinf}
\bv_{\alpha}(\bxi) ~=~ \!\! - \frac{ e^{\text{i}k_p r}\!}{4 \pi(\lambda\!+\!2\mu)r} \exs \bv^\infty_{p_\alpha}(\hat\bxi) \:-\: 
\frac{e^{\text{i}k_s r}\!}{4\pi\mu r} \exs \bv^\infty_{s_\alpha}(\hat\bxi) \:+\: O(r^{-2}) \quad~~ \text{as} \quad r:=|\bxi|\to\infty, 
\vspace*{-1mm}
\eeq
where $\bv^\infty_{p_\alpha}\!\!\parallel\!\hat\bxi$ and $\bv^\infty_{s_\alpha}\!\!\perp\!\hat\bxi$ denote respectively the far-field patterns affiliated with $\bv_{p_\alpha}\!$ and $\bv_{s_\alpha}\!$ in~\eqref{vpvs}, satisfying the integral representations corresponding to~\eqref{GE}-\eqref{KS}, 
\vspace*{-1mm}
\beq\lb{vinfp}
\begin{aligned}
& \bv^\infty_{p_\alpha}(\hat\bxi) ~=~ - \int_{\Gamma_{\!\circ} \cup \Gamma_\alpha} \!\! \llbracket \bv_\alpha \rrbracket(\by) \!\cdot\!  \bSig^\infty_p(\hat\bxi,\by) \!\cdot\! \bn_\alpha(\by) \, \text{d}S_{\by}  ~=~ \\*[0mm] 
- \text{i} k_p \exs \hat\bxi  \int_{\Gamma_{\!\circ} \cup \Gamma_\alpha} \!\! & \Big\lbrace \lambda \,  \llbracket \bv_\alpha \rrbracket \sip \bn_\alpha   + 2\mu \big(\bn_\alpha \sip \hat\bxi \big) \exs \llbracket \bv_\alpha \rrbracket \sip \hat\bxi  \exs \Big\rbrace \, e^{-\text{i}k_p \hat\bxi \cdot \by}  \,\, \text{d}S_{\by}, \quad \alpha \in \lbrace\circ, \kappa\rbrace\nxs, \,\, \hat\bxi\in\Omega, 
\end{aligned} 
\eeq
\vspace*{-5mm}
\beq\lb{vinfs}
\begin{aligned}
& \bv^\infty_{s_\alpha}(\hat\bxi) ~=~ - \int_{\Gamma_{\!\circ} \cup \Gamma_\alpha} \!\! \llbracket \bv_\alpha \rrbracket(\by) \!\cdot\!  \bSig^\infty_s(\hat\bxi,\by) \!\cdot\! \bn_\alpha(\by) \, \text{d}S_{\by}  ~=~ \\*[0mm] 
 -\text{i}k_s \exs \hat\bxi \exs \times \int_{\Gamma_{\!\circ} \cup \Gamma_\alpha} \!\! \Big\lbrace \mu & \big( \llbracket \bv_\alpha \rrbracket \!\times\! \hat\bxi \exs \big)(\bn_\alpha \sip \hat\bxi \exs ) \,+\, \mu \big( \bn_\alpha \!\times \hat\bxi \big) (\llbracket \bv_\alpha \rrbracket \sip \hat\bxi) \Big\rbrace \, e^{-\text{i}k_s \hat\bxi \cdot \by} \,\, \text{d}S_{\by}, \quad \alpha \in \lbrace\circ, \kappa\rbrace\nxs, \,\, \hat\bxi\in\Omega.
\end{aligned} 
\vspace*{-1mm}
\eeq

Here, $\bSig^\infty_p$\! and $\bSig^\infty_s$\! respectively indicate the far-field P- and S- patterns of the elastodynamic fundamental stress tensor $\bSig = \bSig_s \nxs\oplus\nxs \bSig_p$ (see~\cite[Appendix~B]{Fatemeh2017}).

\section{Anatomy of the inverse scattering solution}\lb{Prelim}    

This section introduces key elements of sequential sensing pertinent to  the analysis in \cref{DEI}. 

At every sensing step $t_{\alpha}$, $\alpha \in \lbrace\circ, \kappa\rbrace\nxs$, the domain is excited by a set of plane waves identified with their direction of propagation $\bd \in \Omega$ and polarization amplitudes $\bq=\bq_p\oplus\bq_s$, as in~\eqref{plwa}, and thus-scattered far-field patterns $\bv^\infty_\alpha(\hat\bxi\exs|\exs\bd,\nxs\bq)$ are recorded over a set of observation angles $\hat\bxi \in \Omega$ according to~\eqref{vinf}. In this setting, the far field kernel \textcolor{black}{$\bW_{\!\alpha}^\infty(\bd,\hat\bxi)\in\mathbb{C}^{6\times 6}$} is constructed from far-field data such that 
\vspace*{-1mm}
\beq\lb{w-inf}
\bW_{\!\alpha}^\infty(\bd,\hat\bxi) \sip \bq ~:=~ \bv^\infty_\alpha(\hat\bxi\exs|\exs\bd,\nxs\bq), \qquad \bd,\hat\bxi \in \Omega, \,\, \bq \in \R^3.
\vspace*{-1mm}
\eeq

Given $\bW_{\!\alpha}^\infty\!\!$ at any $t_{\alpha}$, the far field operator $F_{\alpha}: L^2(\OOd)^3 \to L^2(\OOd)^3$ is defined by
\vspace*{-1mm}
\beq\lb{ffo2} 
F_{\alpha}(\bg)\textcolor{black}{(\hat\bxi)} ~=\,  \int_{\OOd} \bW_{\!\alpha}^\infty(\bd,\hat\bxi) \sip \bg(\bd) \,\, \text{d}S_{\bd}, \qquad \alpha \in \lbrace\circ, \kappa\rbrace, \quad\!\! \bg \in L^2(\OO)^3, \quad\!\! \hat\bxi \in \Omega,
\vspace*{-1mm}
\eeq 

Each density function $\bg \in L^2(\OOd)^3$ can be uniquely decomposed as $\bg = \bg_p\!\oplus\bg_s$ where  $\forall\bd\in\Omega$, $\bg_p(\bd)\!\parallel\bd$ and~$\bg_s(\bd)\!\perp\!\bd$. Then, the far field operator maps a density $\bg \in L^2(\OOd)^3$  to the far-field pattern of~$\bv_\alpha \in H^1_{\mathrm{loc}}(\R^3 \backslash \lbrace \Gamma_{\!\circ} \cup \Gamma_\alpha \rbrace)^3\!$ solving~\eqref{GE}-\eqref{KS} when $\bu\ff = \bu_{\bg}$ and where 
\vspace*{-1mm}
\beq\lb{plwa2}
\bu_\bg(\bxi) ~: =~  \int_{\OO} \bg_p(\bd) \exs  e^{\textrm{i} k_p \bd \cdot \bxi} \,\, \text{d}S_{\bd} ~\; \oplus \, \int_{\OO} \bg_s(\bd) \exs  e^{\textrm{i} k_s \bd \cdot \bxi} \,\, \text{d}S_{\bd},  \qquad \bxi \in \R^3, \vspace*{-1mm}
\eeq
denotes a Herglotz wavefield of density $\bg = \bg_p\!\oplus\bg_s$~\cite{Dassios1995}. 
At every sensing step $t = t_\alpha$, we define the Herglotz operator $\mathcal{H}_{\alpha} \colon L^2(\OOd)^3 \rightarrow H^{-1/2}(\Gamma_{\!\circ} \cup \Gamma_\alpha)^3$ that maps the incident polarization densities~$\bg$ of~\eqref{plwa2} to the free-field traction $\bt\ff_\alpha\!$ induced over the scattering interfaces,
\vspace*{-1mm}
\beq\lb{oH}
\mathcal{H}_{\alpha}(\bg) ~:=~ \bn_\alpha\!\cdot\bC\exs\colon\!\nabla\bu_\bg \quad~~ \text{on}\quad \Gamma_{\!\circ} \cup \Gamma_\alpha.
\vspace*{-1mm}
\eeq

With reference to~\eqref{vinfp} and~\eqref{vinfs}, it is then straightforward to show that the adjoint Herglotz operator $\mathcal{H}_\alpha^* \colon \tilde{H}^{1/2}(\Gamma_{\!\circ} \cup \Gamma_\alpha)^3 \rightarrow L^2(\OOd)^3$ takes the form
\vspace*{-1mm}
\beq\lb{Hstar}
\begin{aligned}
& \mathcal{H}_\alpha^*\nxs(\ba)(\hat\bxi) ~=~  - \int_{\Gamma_{\!\circ} \cup \Gamma_\alpha} \!\nxs \ba(\by) \nxs\cdot\nxs  \bSig^\infty(\hat\bxi,\by) \!\cdot\nxs \bn_\alpha(\by) \, \text{d}S_{\by}, \qquad   \bSig^\infty \nxs=\,  \bSig^\infty_s \!\oplus\nxs  \bSig^\infty_p.  
\end{aligned}
\vspace*{-1mm}
\eeq 

Then, each far field operator $F_{\alpha}$ possesses the factorization  
\vspace*{-1mm}
\beq\lb{facts}
F_{\alpha} ~=~ \mathcal{H}_{\alpha}^* \exs T_{\alpha} \exs \mathcal{H}_{\alpha}, \qquad \alpha \in \lbrace\circ, \kappa\rbrace,
\vspace*{-1mm}
\eeq
where the operator $T_{\alpha}\colon H^{-1/2}(\Gamma_{\!\circ} \cup \Gamma_\alpha)^3 \rightarrow \tilde{H}^{1/2}(\Gamma_{\!\circ} \cup \Gamma_\alpha)^3\!$, at $t = t_\alpha$, takes the free field traction $\bt\ff_\alpha\!$ of~\eqref{oH} to the scattered FOD $\ba \colon\!\!\!\!= \llbracket \bv_\alpha \rrbracket$ across $\Gamma_{\!\circ} \cup \Gamma_\alpha$ via the elastic contact laws of~\eqref{GE} and~\eqref{GEi}, 
\vspace*{-2mm}
\beq\lb{T}
T_{\alpha}(\bt\ff_{\alpha})(\bxi) ~:=~ \llbracket \bv_{\alpha}(\bxi) \rrbracket,   \qquad \bxi \in \Gamma_{\!\circ} \cup \Gamma_\alpha, \quad \alpha \in \lbrace\circ, \kappa\rbrace,
\vspace*{-1mm}
\eeq
where~$\bv_\alpha \in H^1_{\mathrm{loc}}(\R^3 \backslash \lbrace \Gamma_{\!\circ} \cup \Gamma_\alpha \rbrace)^3\!$ solves~\eqref{GE}-\eqref{KS} for $\bu\ff = \bu_{\bg}$.

The ensuing analysis requires the following assumption.
\begin{assumption}[illumination prompts scattering]\label{Inject-H}
At any sensing step $t_\alpha$, it is assumed that $\Gamma_{\!\circ} \cup \Gamma_\alpha\!$ and $\omega$ are such that the Herglotz operator $\mathcal{H}_{\alpha} \colon L^2(\OOd)^3 \rightarrow H^{-1/2}(\Gamma_{\!\circ} \cup \Gamma_\alpha)^3$ is injective, and thus, its adjoint $\mathcal{H}_\alpha^* \colon \tilde{H}^{1/2}(\Gamma_{\!\circ} \cup \Gamma_\alpha)^3 \rightarrow L^2(\OOd)^3$ has a dense range. 
\end{assumption}

\cref{Inject-H} is expected to hold true in general for all $\omega\!>\!0$ \textcolor{black}{possibly} excluding a discrete set of values without finite accumulation points. This may be observed by decomposing $\Gamma_{\!\circ} \cup \Gamma_\alpha\!$ into $M_\alpha\!\geqslant\!1$ (possibly disjoint) analytic surfaces~$\Gamma_m\!\subset\Gamma_{\!\circ} \cup \Gamma_\alpha$, $m=1,\ldots M_\alpha$, and identifying their unique analytic continuation $\partial D_m$. In this setting,~\cref{Inject-H} holds according to~\cite[Lemma 5.3]{Fatemeh2017} if $\omega\!>\!0$ is not a ``Neumann'' eigenfrequency of the Navier equation in $D_m$ that satisfies
\vspace*{-2mm}
\beq\lb{uiH}
\begin{aligned}
&\nabla \sip (\bC \colon \! \nabla \pff) \,+\, \rho \exs \omega^2 \pff ~=~ \bzero     \quad  &\textrm{in}~ D_m, \\*[1 mm]
&\bn \sip \bC \colon \! \nxs \nabla  \pff ~=~ \bzero   \quad &\textrm{on}~ \partial D_m.
\end{aligned}
\vspace*{-1mm}
\eeq

Here, $\pff \in H^1(D_m)^3$ indicates the eigen-waveform affiliated to $\omega$. If~$D_m$ is bounded, the real eigenfrequencies of~\eqref{uiH} form a discrete set according to~\cite[Chapter VII, Theorem 1.4]{Kuprad1979}. 

Analogous to~\cite[Lemmas 5.2, 5.6, 5.7]{Fatemeh2017}, it may be shown that at all sensing steps $t_\alpha$:~(a)~$\mathcal{H}_\alpha^*$ is compact and injective, and~(b)~$T_{\alpha}$ is bounded and coercive -- i.e., there exists a constant $c\!>\!0$ independent of~$\bphi$ such that
\vspace{-1 mm}
\beq\lb{co-T0}
|\langle \bphi, \, T_{\alpha} (\bphi) \rangle| \,\,\geqslant\,\, \textrm{c} \nxs \norms{\bphi}_{H^{-\frac{1}{2}}(\Gamma_{\!\circ} \cup \Gamma_\alpha)}^2, \qquad  \forall\bphi\in H^{-1/2}(\Gamma_{\!\circ} \cup \Gamma_\alpha)^3, \vspace{-1.5mm}
\eeq 
and has a continuous inverse. Subsequently, according to~\cite[Lemma 6.6]{Fatemeh2017}, the far-field operator $F_{\alpha}$ is injective, compact and under~\cref{Inject-H} has a dense range.  
Given~\eqref{ffo2}, the operator $F_{\alpha_\sharp} \colon  L^2(\OOd)^3  \rightarrow L^2(\OOd)^3$ is defined by 
\vspace{-1 mm}
\beq\lb{Fs}
F_{\alpha_\sharp}\,\colon \!\!\!=\, \frac{1}{2} |F_{\alpha}+F_{\alpha}^*| \:+\: \frac{1}{2 \textrm{\emph{i}}} (F_{\alpha}-F_{\alpha}^*).
\vspace{-1 mm}
\eeq

Similar to~\cite[Theorem 6.3, Lemma 6.4]{Fatemeh2017}, one may prove the following where the space $\tilde{H}^{1/2}(\Gamma_{\!\circ} \cup \Gamma_\alpha)$ identifies the dual of $H^{-1/2}(\Gamma_{\!\circ} \cup \Gamma_\alpha)$.
\begin{theorem}
Under~\cref{welp} and \cref{Inject-H}, the operator $F_{\alpha_\sharp}$ is positive and has the following factorization
\vspace{-1 mm}
\beq\lb{facts2}
F_{\alpha_\sharp} ~=~ \mathcal{H}_{\alpha}^* \exs T_{\alpha_\sharp} \exs \mathcal{H}_{\alpha}, \qquad \alpha \in \lbrace\circ, \kappa\rbrace,
\vspace{-1 mm}
\eeq
where the middle operator $T_{\alpha_\sharp} \colon H^{-1/2}(\Gamma_{\!\circ} \cup \Gamma_\alpha)^3\rightarrow\tilde{H}^{1/2}(\Gamma_{\!\circ} \cup \Gamma_\alpha)^3$ is 
selfadjoint and positively coercive, i.e., there exists a constant~$\text{\emph{c}}\!>\!0$ independent of~$\bphi$ so that 
\vspace{-1 mm}
\beq\lb{co-T}
\big\langle \bphi, \, T_{\alpha_\sharp}(\bphi) \big\rangle_{H^{-1/2}(\Gamma_{\!\circ} \cup \Gamma_\alpha),\exs\tilde{H}^{1/2}(\Gamma_{\!\circ} \cup \Gamma_\alpha)} \,\,\geqslant\,\, \text{\emph{c}} \nxs \norms{\bphi}_{H^{-\frac{1}{2}}(\Gamma_{\!\circ} \cup \Gamma_\alpha)}^2, \qquad
\forall \,\bphi\in H^{-1/2}(\Gamma_{\!\circ} \cup \Gamma_\alpha)^3.
\vspace{-1 mm}
\eeq

 Moreover, the range of $\mathcal{H}_\alpha^*$ coincides with that of $F_{\alpha_\sharp}^{1/2}$.
\end{theorem}

{\bfseries\slshape Philosophy of the sampling-based data inversion.}
\textcolor{black}{With reference to~\cref{comp_ps}, let us define a search volume $\mathcal{B} \subset \R^3$ and a set of trial scatterers $L(\bx_\circ, \bR) \subset \mathcal{B}$ such that for every pair $(\bx_\circ,\bR)$, $L\colon\!\!=\bx_\circ\!+\bR{\sf L}$ specifies a smooth arbitrary-shaped dislocation~$\sf L$~at $\bx_\circ \subset \mathcal{B}$ whose orientation is identified by a unitary rotation matrix $\bR\!\in\!U(3)$. In this setting, the far-field pattern $\bPhi_L^\infty \colon \tilde{H}^{1/2}(L) \rightarrow L^2(\OOd)^3$ induced by $L(\bx_\circ, \bR)$ as a sole scatterer in $\R^3$ endowed with an admissible FOD profile $\ba\!\in\!\tilde{H}^{1/2}(L)$ is given by} 
\vspace{-1 mm}
\beq\lb{Phi-inf}
\begin{aligned}
\bPhi_{L}^\infty(\ba)(\hat\bxi) ~=~ -  \Big(  \, & \text{\emph{i}} k_p \,  \hat\bxi \exs \int_L \, \Big\lbrace \lambda \exs (\ba \sip \bn) \,+\, 2\mu \exs (\bn \sip \hat\bxi) ( \ba \sip \hat\bxi)  \Big\rbrace  \, e^{-\text{\emph{i}}k_p \hat\bxi \cdot \by} \,\, \text{d}S_{\by}  \\*[1 mm]
 & \textcolor{black}{\oplus}~ \text{\emph{i}} k_s \,  \hat\bxi\times\! \int_L \Big\lbrace   \mu \exs(\ba \times \hat\bxi)(\bn\sip\hat\bxi) \,+\, \mu \exs  (\bn \times \hat\bxi) (\ba \sip \hat\bxi)  \Big\rbrace \, e^{-\textrm{\emph{i}} k_s \hat\bxi \cdot \by} \,\, \text{d}S_{\by} \Big),
\end{aligned}
\vspace{-1 mm}
\eeq
where $\bn(\by)$ signifies the unit normal at $\by \in L$, and $\hat\bxi \in \Omega$ is the observation direction. In light of~\eqref{Phi-inf}, one may generate a library of  scattering signatures affiliated with a grid of trial pairs $(\bx_\circ,\bR)$ sampling $\mathcal{B}\!\times\! U(3)$. 

~The underpinning concept of sampling methods~\cite{audi2015, Fatemeh2017} is segment-wise reconstruction of $\Gamma_{\!\circ} \cup \Gamma_\alpha$ through a careful implementation of synthetic wavefront shaping at every sensing step $t_\alpha$, aiming to distill the scattering signature of domain's internal structure segment by segment from waveform data. In this vein, the library of far field patterns $\bPhi_{L}^\infty$, generated on the basis of trial dislocations $L(\bx_\circ, \bR)$, is deployed to probe the range of the operator $F_\alpha$ (or $F_{\alpha_\sharp}^{1/2}$ in the factorization method) by solving 
\vspace{-1 mm}
\beq\lb{FF}
F_\alpha\exs \bg ~\simeq~ \bPhi_L^\infty, \qquad (\mbox{or } F_{\alpha_\sharp}^{1/2} \bg ~=~ \bPhi_L^\infty \mbox{ for the factorization method})
\vspace{-1 mm}
\eeq
for the illumination densities $\bg(\bx_\circ, \bR)=\bg_p\oplus\bg_s$.
In this setting, the main theorems underlying sampling methods (e.g., the factorization method and generalized linear sampling method) rigorously furnish the distinct behavior of solution $\bg(\bx_\circ, \bR)$ in the vicinity of hidden scatterers, giving rise to a suit of imaging criteria to reconstruct $\Gamma_{\!\circ} \cup \Gamma_\alpha$. We refer to the following section for the indicator related to the generalized linear sampling method. For the factorization method, the equation is solvable if and only if $L \subset \Gamma_{\!\circ} \cup \Gamma_\alpha$.

\begin{rem}[finite domains]\lb{FD}
It should be mentioned that one may also rigorously define parallel operators pertinent to finite backgrounds that carry similar properties mentioned in this section, see e.g.,~\cite{nguy2019}. Thus, the ensuing developments directly lend themselves to both finite and infinite domains, as demonstrated by the numerical experiments in \cref{IR}.   
\end{rem}

\section{Differential evolution indicators}\lb{DEI}  
 \vspace{-2mm}

\begin{figure}[h!]
\vspace*{-2mm}
\begin{center}
\includegraphics[width=1\linewidth]{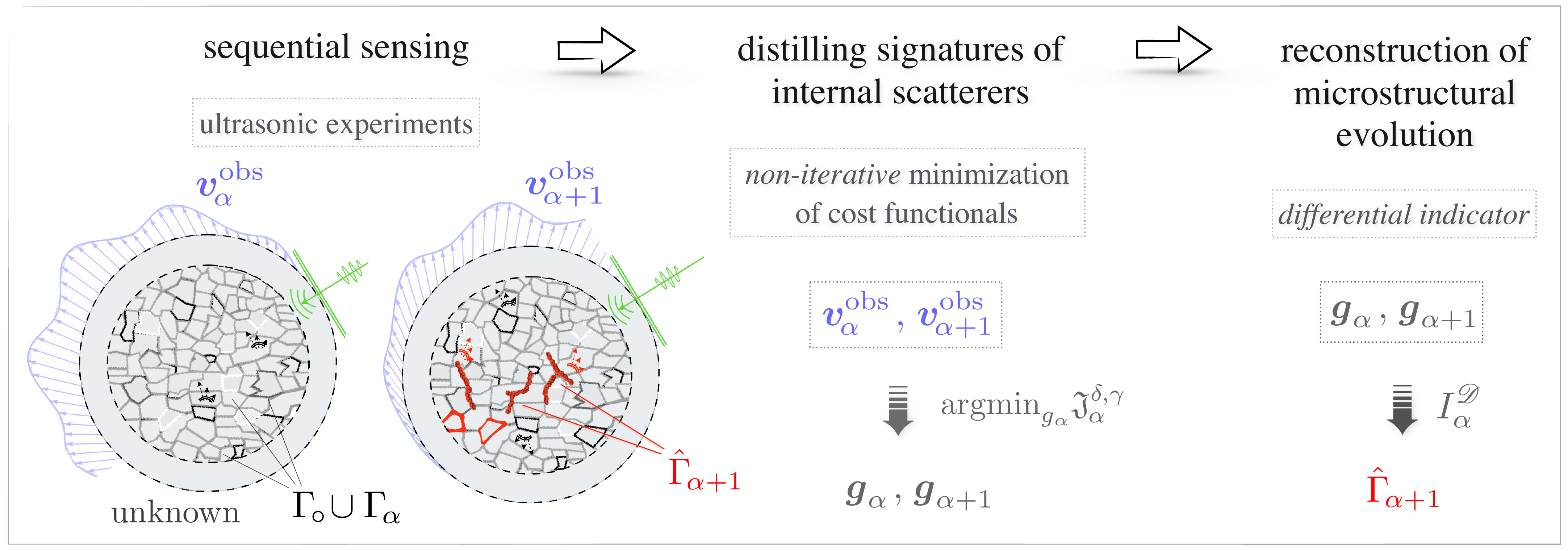}
\end{center} \vspace*{-3.2mm}
\caption{\small{Three-tier approach to differential tomography of microstructural evolution in highly scattering backgrounds.}}
\label{inversion}\vspace*{-2.5mm}
\end{figure}

As mentioned earlier, existing sampling approaches to wavefrom tomography mostly require an accurate characterization of the background for their successful performance. This section aims to relax this requirement by introducing a \emph{three-tier platform} for targeted reconstruction of evolution in elastic backgrounds of a-priori unknown structure. As illustrated in~\cref{inversion}, in this framework:~{\bf (1)}~the domain is sequentially subject to illumination and sensing in an arbitrary configuration;~{\bf (2)}~the resulting sequence of sensory data~$(\bv_{\alpha}^{\textrm{obs}}, \bv_{\alpha+1}^{\textrm{obs}})$ are deployed to \emph{non-iteratively} compute the associated source densities i.e., synthetic wavefronts $(\bg_\alpha, \bg_{\alpha+1})$ to distill segment-wise signatures of the domain's internal structure from scattered field measurements; and~{\bf (3)}~thus-obtained densities are then used to selectively reconstruct the support of interfaces $\hat{\Gamma}_{\alpha+1} \colon\!\!\!\!=\! {\Gamma}_{\alpha+1} \backslash {\Gamma}_{\alpha}$ born or evolved between $t \in [t_\alpha, t_{\alpha+1}]$ (or any pairs of sensing steps) within a network of pre-existing scatterers ${\Gamma}_{\!\alpha} \cup {\Gamma}_{\!\circ}$.   

Construction of an evolution indicator is rooted in minimizing sequences to the cost functional 
\vspace*{-2mm}
\beq\lb{GCf}
\mathfrak{J}_\alpha^\gamma(\bPhi_L^\infty;\,\bg) ~:=~ \norms{\nxs F_\alpha\exs \bg\,-\,\bPhi_L^\infty \nxs}_{L^2(\Omega)}^2 +~ \gamma\, \big(\bg, F_{\alpha_\sharp}\bg\big), \qquad \bg\in L^2(\OOd)^3, \,\,\, \gamma>0.
\vspace*{-1mm}
\eeq
On denoting $G_\alpha = \mathcal{H}^*_\alpha \exs T_\alpha$, we also define 
$$ \mathfrak{L}_\alpha^\gamma(\bPhi_L^\infty; \bpsi) = \norms{ G_\alpha\bpsi \,-\,\bPhi_L^\infty }_{L^2(\Omega)}^2 +~ \gamma\, \big(\bpsi, T_{\alpha_\sharp}\bpsi\big), \qquad \bpsi \in H^{-\frac{1}{2}}(\Gamma_{\!\circ} \cup \Gamma_\alpha)^3, \,\,\, \gamma>0,$$
where $$\mathfrak{L}_\alpha^\gamma(\bPhi_L^\infty; \mathcal{H}_{\alpha} \bg) = \mathfrak{J}_\alpha^\gamma(\bPhi_L^\infty;\bg).$$

In what follows, the strong convergence of germane minimizing sequences is established by way of the strong convexity of $\mathfrak{L}_\alpha^\gamma(\bPhi_L^\infty; \cdot) $ on ${H^{-\frac{1}{2}}(\Gamma_{\!\circ} \cup \Gamma_\alpha)}$. This approach is slightly different from the related arguments in~\cite{audiDLSM}. 

\begin{theorem} \lb{GLSM1}
Consider the minimizing sequence $\bg^\gamma \in L^2(\OOd)^3$ for $\mathfrak{J}_\alpha^\gamma$ such that   
\vspace{-1mm}
\beq\lb{mseq1}
\mathfrak{J}_\alpha^\gamma(\bPhi_L^\infty;\bg^\gamma) \,\leqslant \, \mathfrak{j}_\alpha^\gamma(\bPhi_L^\infty) + \eta(\gamma), \qquad \gamma > 0, \vspace{-1mm}
\eeq
where $\eta(\gamma)/\gamma \rightarrow 0$ as $\gamma \rightarrow 0$ and
\vspace{-1mm}
\[
\mathfrak{j}_\alpha^\gamma(\bPhi_L^\infty) ~\colon\!\!\!=~\!\! \inf\limits_{\bg \in L^2(\OOd)^3} \! \mathfrak{J}_\alpha^\gamma(\bPhi_L^\infty;\bg). 
\vspace{-4mm}
\]

Then,  
  \vspace{-1mm}
\beq\lb{statG}
\begin{aligned}
& \bPhi_L^\infty \in Range(G_\alpha) ~\Rightarrow~ \lim\limits_{\gamma \rightarrow 0} \big(\exs \bg^\gamma, F_{\alpha_\sharp}\bg^\gamma \exs \big) = \lim\limits_{\gamma \rightarrow 0} \big(\exs \mathcal{H}_\alpha \exs \bg^\gamma, T_{\alpha_\sharp} \mathcal{H}_\alpha \exs \bg^\gamma \exs \big) < \infty, \\
&\bPhi_L^\infty \notin Range(G_\alpha) ~\Rightarrow~ \liminf\limits_{\gamma \rightarrow 0} \big(\exs \bg^\gamma, F_{\alpha_\sharp}\bg^\gamma \exs \big) = \liminf\limits_{\gamma \rightarrow 0} \big(\exs \mathcal{H}_\alpha \exs \bg^\gamma, T_{\alpha_\sharp}\mathcal{H}_\alpha \exs \bg^\gamma \exs \big) = \infty.
 \end{aligned}
  \vspace{-2mm}
\eeq

 Moreover, when $G_\alpha \exs \boldsymbol{\psi} =  \bPhi_L^\infty$, the sequence $\mathcal{H}_\alpha \bg^\gamma$ strongly converges to $\boldsymbol{\psi} \in {H^{-\frac{1}{2}}(\Gamma_{\!\circ} \cup \Gamma_\alpha)}$ as $\gamma \rightarrow 0$. 
\end{theorem}

\begin{proof}
The limits of~\eqref{statG} are directly drawn from~\cite[Theorem 6.7]{Fatemeh2017}. In the case where  $G_\alpha \exs \boldsymbol{\psi} =  \bPhi_L^\infty$ for some $ \boldsymbol{\psi} \in {H^{-\frac{1}{2}}(\Gamma_{\!\circ} \cup \Gamma_\alpha)}$,
the limiting behavior of $\mathcal{H}_\alpha \exs \bg^\gamma$ to $\boldsymbol{\psi}$ may be  analyzed by using the strong convexity  of $\mathfrak{L}_\alpha^\gamma(\bPhi_L^\infty; \cdot) $. More specifically, using that $\mathfrak{J}_\alpha^\gamma(\bPhi_L^\infty;\bg^\gamma) = \mathfrak{L}_\alpha^\gamma(\bPhi_L^\infty;\mathcal{H}_{\alpha}\bg^\gamma)$ and that  $\mathcal{H}_{\alpha}$ has dense range, we have 
\vspace{-1mm}
\beq \lb{StConv}
\begin{aligned}
& \mathfrak{J}_\alpha^\gamma(\bPhi_L^\infty;\bg^\gamma) \,-\, \eta(\gamma) \,+\, \gamma \exs {\textstyle \frac12} \norms{\nxs\mathcal{H}_{\alpha}\bg^\gamma-\boldsymbol{\psi}\nxs}^2 \,\,\,\leqslant\,  \mathfrak{j}_\alpha^\gamma(\bPhi_L^\infty) \,+\, \gamma \exs {\textstyle \frac12} \norms{\nxs\mathcal{H}_{\alpha}\bg^\gamma-\boldsymbol{\psi}\nxs}^2 \,\,\,\leqslant\, \\*[0.5mm]
& \leqslant\, \mathfrak{L}_\alpha^\gamma\big(\bPhi_L^\infty;  {\textstyle \frac12}\mathcal{H}_{\alpha}\bg^\gamma + {\textstyle \frac12} \boldsymbol{\psi}\big) + \gamma \exs {\textstyle \frac12} \norms{\nxs\mathcal{H}_{\alpha}\bg^\gamma-\boldsymbol{\psi}\nxs}^2 \,\, \leqslant \, \max \lbrace \mathfrak{J}_\alpha^\gamma\big(\bPhi_L^\infty;\,\bg^\gamma\big), \, \mathfrak{L}_\alpha^\gamma(\bPhi_L^\infty;\boldsymbol{\psi}) \rbrace.
\end{aligned}
\vspace{-1mm}
\eeq 

Then, in light of 
\vspace{-1mm}
\beq\lb{JJTT}
\mathfrak{L}_\alpha^\gamma(\bPhi_L^\infty;\boldsymbol{\psi}) - \, \mathfrak{J}_\alpha^\gamma(\bPhi_L^\infty;\bg^\gamma) \,\,\leqslant\, \gamma \exs \big[ \langle \boldsymbol{\psi}, \, T_{\alpha_\sharp}\boldsymbol{\psi} \rangle - \langle \mathcal{H}_\alpha \exs \bg^\gamma, \, T_{\alpha_\sharp}\mathcal{H}_\alpha \exs \bg^\gamma \rangle \big], 
\vspace{-1mm}
\eeq
for $\mathfrak{L}_\alpha^\gamma(\bPhi_L^\infty;\boldsymbol{\psi}) > \, \mathfrak{J}_\alpha^\gamma(\bPhi_L^\infty;\bg^\gamma)$, observe from~\eqref{StConv} that
\vspace{-1mm}
\[
\limsup\limits_{\gamma \rightarrow 0} \norms{\nxs\mathcal{H}_{\alpha}\bg^\gamma-\boldsymbol{\psi}\nxs}^2 = 0,
\vspace{-2mm}
\] 
which proves the strong convergence of $\mathcal{H}_{\alpha}\bg^\gamma$ to $\boldsymbol{\psi}$ as $\gamma \rightarrow 0$.       
\end{proof}

{\bfseries\slshape Noisy data.} Consider the perturbed operators $F^\delta_\alpha\nxs$ and $F^\delta_{\alpha_\sharp}\nxs$, 
\vspace{-1mm}
\beq\lb{Ns-op}
\norms{\nxs F^\delta_{\alpha} - F_{\alpha} \nxs} \,\, \leqslant \, \delta, \qquad \norms{\nxs F^\delta_{\alpha_\sharp} - F_{\alpha_\sharp} \nxs} \,\, \leqslant \, \delta,   \vspace{-1mm}
\eeq
where $\delta\!>\!0$  is a measure of noise in data and the self-adjoint operator $F^\delta_{\alpha_\sharp} \colon  L^2(\OOd)^3  \rightarrow L^2(\OOd)^3$ is drawn from $F_{\alpha}^\delta$ via 
\vspace{-2mm}
\beq\lb{Fsd}
F^\delta_{\alpha_\sharp}\,\colon \!\!\!=\, \frac{1}{2} \big{|}F^\delta_{\alpha}+{F^\delta_{\alpha}}^*\big{|} \:+\: \big{|}\frac{1}{2 \textrm{\emph{i}}} (F^\delta_{\alpha}-{F^\delta_{\alpha}}^*)\big{|}. 
\vspace{-5 mm}
\eeq
\begin{assumption}\label{ND}
$\forall t_\alpha$, $F^\delta_\alpha\nxs$ and $F^\delta_{\alpha_\sharp}\nxs$ are compact.
\end{assumption}

{
\begin{theorem}[noisy data] \lb{GLSM2} 
For $\bg \in L^2(\OOd)^3, \,\, \gamma\!>\!0$, consider the cost functional  
\vspace{-1mm}
\beq\lb{GCfn}
\mathfrak{J}_\alpha^{\delta, \gamma}(\bPhi_L^\infty;\bg) \,\,:=\,\,\,\exs  \norms{\nxs F_{\alpha}^\delta\exs\bg\,-\,\bPhi_L^\infty \nxs}_{L^2(\Omega)}^2 \!\,+\,\,\exs\gamma\, \big(\bg, F^\delta_{\alpha_\sharp}\bg\big) \,+\, \gamma^{1-\chi}\exs\delta  \norms{\nxs \bg \nxs}^2_{L^2(\Omega)} ,
\vspace{-1mm}
\eeq
where $\chi \in \, ]0, 1[$, and $\mathfrak{J}_\alpha^{\delta, \gamma}(\bPhi_L^\infty;\bg)$ admits the minimizer  
\vspace{-1mm}
\beq\lb{gL}
\bg^{\delta,\gamma}_L ~\!=~ \text{\emph{arg}}\!\!\!\!\!\min_{\bg \in L^2(\OOd)^3}  \mathfrak{J}_\alpha^{\delta, \gamma}(\bPhi_L^\infty;\bg),   \vspace{-1mm}
\vspace{-1.5mm}
\eeq 
satisfying
  \vspace{-1mm}
\beq\lb{limlim}
\lim\limits_{\gamma \rightarrow 0} \limsup\limits_{\delta \rightarrow 0} \exs  \mathfrak{J}_\alpha^{\delta, \gamma}(\bPhi_L^\infty;\bg^{\delta,\gamma}_L) ~\!=\!~ 0.   \vspace{-1.5mm}
\eeq 
In this setting,
  \vspace{-1mm}
\beq\lb{statG2}
\begin{aligned}
& \bPhi_L^\infty \in Range(G_\alpha) ~\Rightarrow~ \limsup\limits_{\gamma \rightarrow 0}\limsup\limits_{\delta \rightarrow 0}\big( (\exs \bg^{\delta,\gamma}_L, F^\delta_{\alpha_\sharp}\, \bg^{\delta,\gamma}_L) \,+\, \delta \gamma^{-\chi} \! \norms{ \bg^{\delta,\gamma}_L }^2  \!\big) \,<\, \infty, \\
&\bPhi_L^\infty \notin Range(G_\alpha) ~\Rightarrow~ \liminf\limits_{\gamma \rightarrow 0}\liminf\limits_{\delta \rightarrow 0}\big( (\exs \bg^{\delta,\gamma}_L, F^\delta_{\alpha_\sharp}\, \bg^{\delta,\gamma}_L) \,+\, \delta \gamma^{-\chi} \! \norms{ \bg^{\delta,\gamma}_L }^2 \! \big) \,=\, \infty.
 \end{aligned}
  \vspace{-1mm}
\eeq
In addition, when $G_\alpha \boldsymbol{\psi} =  \bPhi_L^\infty$, there holds
\vspace{-2mm}
\beq\lb{lmt}
\limsup\limits_{\gamma \rightarrow 0} \limsup\limits_{\delta \rightarrow 0} \exs \delta \! \norms{ \bg^{\delta,\gamma}_L }^2 ~\!\!=\!~ 0. \vspace{-2mm}
\eeq
\vspace{-2mm}
Also, there exists $\delta_{\circ}(\gamma)$ such that $\forall \delta(\gamma) \leqslant \delta_{\circ}(\gamma)$, $\mathcal{H}_{\alpha}\bg^{\delta(\gamma),\gamma}_L\!$ converges strongly to $\boldsymbol{\psi}$ as $\gamma \rightarrow 0$. \end{theorem} 
}
\vspace{0mm}
\begin{proof}
{The limiting behavior of $\mathfrak{J}_\alpha^{\delta, \gamma}$ in~\eqref{limlim}, and limits of the penalty term in~\eqref{statG2} are established in \cite{audiDLSM}. Moreover, given~\eqref{statG2} for the case where $\bPhi_L^\infty \in Range(G_\alpha)$,~\eqref{lmt} is self-evident.} This relation along with the strong convergence result, when $G_\alpha \boldsymbol{\psi} =  \bPhi_L^\infty$, constitutes the foundation of differential imaging with noisy operators and may be observed as the following. Define 
  \vspace{-1mm}
\beq\lb{d(g)}
\begin{aligned}
& {D}_\alpha^{\delta, \gamma}(\bPhi_L^\infty;\bg) \, \colon\!\!\!\nxs= \, \delta^2 \!\nxs \norms{\nxs \bg \nxs}^2 \nxs+\, 2 \delta \! \norms{\nxs F_{\alpha} \exs\bg\,-\,\bPhi_L^\infty \nxs} \norms{\nxs \bg \nxs} \nxs+\, \delta \gamma \! \norms{\nxs \bg \nxs}^2, \\*[0.5mm]
& \tilde{D}_\alpha^{\delta, \gamma}(\bPhi_L^\infty;\bg) \, \colon\!\!\!\nxs= \, {D}_\alpha^{\delta, \gamma}(\bPhi_L^\infty;\bg)  \nxs+\, \delta \gamma^{1-\chi} \! \norms{\nxs \bg \nxs}
^2.
\end{aligned}
  \vspace{-1mm}
\eeq

Then, observe that
  \vspace{-2mm}
\beq\lb{diffJJ}
\begin{aligned}
&[\exs\mathfrak{J}_\alpha^{\delta, \gamma} -\, \mathfrak{J}_\alpha^{\gamma}\exs](\bPhi_L^\infty;\bg) \, \leqslant \, \tilde{D}_\alpha^{\delta, \gamma}(\bPhi_L^\infty;\bg), \\*[0.5mm]
&[\mathfrak{J}_\alpha^{\gamma} - \exs\mathfrak{J}_\alpha^{\delta, \gamma} \exs](\bPhi_L^\infty;\bg) \, \leqslant \, {D}_\alpha^{\delta, \gamma}(\bPhi_L^\infty;\bg).
\end{aligned}
  \vspace{-1mm}
\eeq
Consider $\eta(\gamma)$ and $\bg^{\gamma}$ as in~\cref{GLSM1},~\eqref{mseq1},~\eqref{gL} and~\eqref{diffJJ}, then observe that
  \vspace{-1mm}
\beq\lb{StrgCon2}
\begin{aligned}
& \mathfrak{J}_\alpha^{\delta, \gamma}(\bPhi_L^\infty;\bg_L^{\delta, \gamma}) \,-\, \tilde{D}_\alpha^{\delta,\gamma}(\bPhi_L^\infty, \bg^\gamma)\,-\, \eta(\gamma) \,+\, \gamma \exs {\textstyle \frac12} \norms{\nxs\mathcal{H}_{\alpha}\bg_L^{\delta, \gamma}-\boldsymbol{\psi}\nxs}^2 \,\,\,\leqslant\, \\*[0.5mm] 
& \mathfrak{J}_\alpha^{\delta, \gamma}(\bPhi_L^\infty;\bg^{\gamma}) \,-\, \tilde{D}_\alpha^{\delta,\gamma}(\bPhi_L^\infty, \bg^\gamma) \,-\, \eta(\gamma) \,+\, \gamma \exs {\textstyle \frac12} \norms{\nxs\mathcal{H}_{\alpha}\bg_L^{\delta, \gamma}-\boldsymbol{\psi}\nxs}^2 \,\,\,\leqslant\, \\*[0.5mm] 
& \mathfrak{j}_\alpha^\gamma(\bPhi_L^\infty) \,+\, \gamma \exs {\textstyle \frac12}  \norms{\nxs\mathcal{H}_{\alpha}\bg_L^{\delta, \gamma}-\boldsymbol{\psi}\nxs}^2 \,\,\,\leqslant\,   \max \lbrace \mathfrak{J}_\alpha^\gamma\big(\bPhi_L^\infty;\,\bg_L^{\delta, \gamma}\big), \, \mathfrak{L}_\alpha^\gamma\big(\bPhi_L^\infty;\,\boldsymbol{\psi}) \rbrace.
\end{aligned}
  \vspace{-1mm}
\eeq
which in light of \eqref{diffJJ} results in  
$$
\gamma \exs {\textstyle \frac12} \! \norms{\nxs\mathcal{H}_{\alpha}\bg_L^{\delta, \gamma}-\boldsymbol{\psi}\nxs}^2 \exs \leqslant \!  \max \lbrace  \mathfrak{L}_\alpha^\gamma\big(\bPhi_L^\infty;\boldsymbol{\psi})  - \mathfrak{J}_\alpha^{\gamma}(\bPhi_L^\infty;\bg_L^{\delta, \gamma}), 0 \rbrace + {D}_\alpha^{\delta,\gamma}(\bPhi_L^\infty; \bg_L^{\delta, \gamma})+  \tilde{D}_\alpha^{\delta,\gamma}(\bPhi_L^\infty; \bg^\gamma) + \eta(\gamma),
$$
{By definition of ${D}_\alpha^{\delta_{\circ}(\gamma),\gamma}(\bPhi_L^\infty, \bg)$ in~\eqref{d(g)}, it is evident that one may choose a sequence $\delta_{\circ}(\gamma)$ such that $\forall \delta(\gamma) \leqslant \delta_{\circ}(\gamma)$,
  \vspace{-4mm}
\beq\nonumber
\lim\limits_{\gamma \rightarrow 0} \exs \frac{\tilde{D}_\alpha^{\delta(\gamma),\gamma}(\bPhi_L^\infty, \bg^\gamma)}{\gamma} \,=\, 0.
  \vspace{-1mm}
\eeq
In view of \eqref{statG2} and \eqref{lmt}, one may also find $\delta_{\circ}(\gamma)$ such that $\forall \delta(\gamma) \leqslant \delta_{\circ}(\gamma)$,
  \vspace{-1mm}
\beq\nonumber
\lim\limits_{\gamma \rightarrow 0} \exs \frac{{D}_\alpha^{\delta(\gamma),\gamma}(\bPhi_L^\infty, \bg_L^{\delta(\gamma),\gamma})}{\gamma} \,=\, 0.
  \vspace{-1mm}
\eeq
Consequently, from~\eqref{JJTT}, one may conclude that 
  \vspace{-1mm}
\beq\nonumber
\limsup\limits_{\gamma \rightarrow 0} \norms{\nxs\mathcal{H}_{\alpha}\bg_L^{\delta(\gamma), \gamma}-\boldsymbol{\psi}\nxs}^2 \,=\, 0, \qquad \forall \delta(\gamma) \leqslant \delta_{\circ}(\gamma),
\vspace{-4mm}
\eeq
verifying the strong convergence of $\mathcal{H}_{\alpha}\bg_L^{\delta, \gamma}$ to $\boldsymbol{\psi}$ as $\gamma \rightarrow 0$.}
\end{proof}

\vspace{-1mm}
{\bfseries\slshape Single-step GLSM imaging criteria.}~For future reference, it should be mentioned that~\cref{GLSM2} forms the foundation of GLSM imaging indicator~\cite{Fatemeh2017}, 
\vspace{-0.5mm}
\beq\lb{GLSMgs}
I^{{\mathcal{G}}_\sharp}(L) \,\, = \,\, \dfrac{1}{\sqrt{\norms{\!(F^\delta_{\alpha_\sharp})^{\frac{1}{2}} \exs \bg_L^{\delta, \gamma} \nxs}^2 \exs+\,\, \delta \! \norms{ \bg_L^{\delta, \gamma} \!}^2}}, \vspace{-0.5mm}
\eeq      
constructed on the basis of scattered field data $F^\delta_{\alpha_\sharp}$ captured at a single time step $t_\alpha$. $I^{{\mathcal{G}}_\sharp}(L)$ attains its highest values when the trial dislocation $L(\bx_\circ,\bR)$ meets the support of hidden scatterers $\Gamma_{\!\circ} \cup \Gamma_\alpha$.

\begin{rem}
The GLSM indicator $I^{{\mathcal{G}}_\sharp}$ is primarily designed for imaging in elastic backgrounds whose topology and material properties are precisely identified~\cite{audi2015, Fatemeh2017, de2018}. Such rigorous knowledge of the background is not achievable in many practical situations particularly at micro- and meso- scales. Furthermore, in a fully characterized background domain, it is shown in~\cref{IR} that $I^{{\mathcal{G}}_\sharp}$ loses its resolution in presence of multiple closely spaced scatterers whose pairwise distances are of the order of a fraction of the illuminating wavelength. 
\end{rem}

\vspace{-1mm}
{\bfseries\slshape Invariants of scattering solution.}~Requiring an exact knowledge of background may be relaxed by taking advantage of (a) unique attributes of the cost functionals $\mathfrak{J}_\alpha^{\gamma}$ (\emph{resp.}~$\mathfrak{J}_\alpha^{\delta, \gamma}$) introduced in~\cref{GLSM1} (\emph{resp.}~\cref{GLSM2}) -- namely, their convex nature and robustness against noise, and (b) newly established strong convergence of the proposed minimizing sequence $\bg^{\delta,\gamma}_{L}$ to a unique minimizer when $\bPhi_L^\infty \in Range(G_\alpha)$. This claim is further motivated by the following~\cref{TR1}, where the relation between any pairs of synthetic wavefronts $(\bg_\alpha, \bg_{\alpha+1})$, computed in {\emph{Tier 2}}, is established in terms of their affiliated FODs i.e.,~fracture opening displacement profiles $(\dbva, \dbvap)$. FODs are directly linked to the penalty terms in $(\mathfrak{J}^{\delta, \gamma}_\alpha, \mathfrak{J}^{\delta, \gamma}_{\alpha+1})$, and thus imaging indicators of the sampling type e.g.,~$I^{{\mathcal{G}}_\sharp}$ in~\eqref{GLSMgs}. Based on these developments,~\cref{INV1} and~\cref{INV2} introduce a new class of functionals that remain systematically invariant with respect to the stationary scatterers ${\Gamma}_{\!\alpha} \cup {\Gamma}_{\!\circ}$ between any pair of distinct experimental campaigns $t \in [t_\alpha, t_{\alpha+1}]$. This leads to the differential evolution indicators~$I^{\mathcal{D}}_{\alpha}, \hat{I}^{\mathcal{D}}_{\alpha}$ (\emph{resp.}~$I^{\mathcal{D}^\delta}_{\alpha,\delta}, \hat{I}^{\mathcal{D}^\delta}_{\alpha,\delta}$) in~\eqref{EIF} (\emph{resp.}~\eqref{EIFn}) enabling selective reconstruction of evolution $\hat{\Gamma}_{\alpha+1} = {\Gamma}_{\alpha+1} \backslash {\Gamma}_{\alpha}$ within the interval $[t_\alpha, t_{\alpha+1}]$ in a complex background without the need to reconstruct the entire domain i.e., $\Gamma_{\!\circ} \cup \Gamma_\alpha$ across pertinent scales, which may be practically insurmountable.

\begin{theorem}\lb{TR1} 

Given $(\bv_{\alpha}^{\infty}, \bv_{\alpha+1}^{\infty})$, consider sampling the search volume $\bx_\circ \in \mathcal{B}$ by a set of trial dislocations $L(\bx_\circ, \bR)$ endowed with an admissible FOD~$\,\ba(\bxi)\in\tilde{H}^{1/2}(L)^3$. The resulting source densities $(\bg_\alpha, \bg_{\alpha+1})(L;\delta,\gamma)$ minimizing $(\mathfrak{J}^{\delta, \gamma}_\alpha, \mathfrak{J}^{\delta, \gamma}_{\alpha+1})$ are deployed to identify the affiliated Herglotz incidents $(\bu_{\bg_\alpha},\bu_{\bg_{\alpha+1}})$ in~\eqref{plwa2}, and thereof, the scattered FOD profiles $\,\dbva(\bxi)\in\tilde{H}^{1/2}(\Gamma_{\!\alpha} \cup {\Gamma}_{\!\circ})^3$ and $\,\dbvap(\bxi)\in\tilde{H}^{1/2}({\Gamma}_{\alpha+1} \cup {\Gamma}_{\!\circ})^3$ according to~\eqref{GEi}. Then, under~\cref{welp} and~\cref{Inject-H},
\vspace{1mm}
\begin{itemize}
\item{} If $L\subset\Gamma_{\!\alpha} \cup {\Gamma}_{\!\circ}$ then $\,\dbva = \,\dbvap$ over $\Gamma_{\!\alpha} \cup {\Gamma}_{\!\circ}$.
\item{} If $L\subset\hat{\Gamma}_{\!\alpha+1}\!$ then $\,\dbva \neq \,\dbvap = \bzero$ over $\Gamma_{\!\alpha} \cup {\Gamma}_{\!\circ}$.
\end{itemize}
\end{theorem}
\begin{proof}
Consider the following:
\vspace{1mm}
\begin{itemize}
\item If~$L \subset \Gamma_{\!\alpha} \cup {\Gamma}_{\!\circ}$, then $\tilde{H}^{1/2}(L)^3 \subset \tilde{H}^{1/2}(\Gamma_{\!\alpha} \cup {\Gamma}_{\!\circ})^3$ (resp.~$\tilde{H}^{1/2}(L)^3 \subset \tilde{H}^{1/2}(\Gamma_{\!\alpha+1}\! \cup {\Gamma}_{\!\circ})^3$). By extending the domain of $\ba\in\tilde{H}^{1/2}(L)^3$ from $L$ to $\Gamma_{\!\alpha} \cup {\Gamma}_{\!\circ}$ (resp.~$\Gamma_{\!\alpha+1} \cup {\Gamma}_{\!\circ}$) through zero padding, one immediately obtains $\bPhi_L^\infty  \in Range(\mathcal{H}_\alpha^*)$ (resp.~$\bPhi_L^\infty  \in Range(\mathcal{H}_{\alpha+1}^*)$) thanks to~\eqref{Hstar} and~\eqref{Phi-inf}. As a result, $\mathcal{H}_{\alpha}^* \dbva = \bPhi_L^\infty$ (resp.~$\mathcal{H}_{\alpha+1}^* \dbvap = \bPhi_L^\infty$) possesses a uniques solution such that $\dbva = \ba$ over $L$ and $\dbva = \bzero$ on $\Gamma_{\!\alpha} \cup {\Gamma}_{\!\circ} \backslash L$ (resp.~$\dbvap = \ba$ over $L$ and $\dbvap= \bzero$ on $\Gamma_{\alpha+1} \cup {\Gamma}_{\!\circ} \backslash L$).  Therefore, $\,\dbva = \,\dbvap$ over $\Gamma_{\!\alpha}  \cup {\Gamma}_{\!\circ}$.

\item If~$L\subset\hat{\Gamma}_{\!\alpha+1}$, then $\bPhi_L^\infty\!\in\!Range(\mathcal{H}_{\alpha+1}^*)$ while $\bPhi_L^\infty\!\not\in\!Range(\mathcal{H}_{\alpha}^*)$. In this case, according to the above argument $\mathcal{H}_{\alpha+1}^* \dbvap = \bPhi_L^\infty$ has a unique solution, and $\dbvap = \bzero$ on $\Gamma_{\!\alpha} \cup {\Gamma}_{\!\circ}$ as it falls in the zero-padded region (recall that $L \not\subset\Gamma_{\!\alpha} \cup {\Gamma}_{\!\circ}$). On the other hand, the norm of any approximate solution to $\mathcal{H}_{\alpha}^* \dbva = \bPhi_L^\infty$ become unbounded ($\norms{\!\dbva\!}_{\tilde{H}^{1/2}(\Gamma_{\!\alpha} \cup {\Gamma}_{\!\circ})^3} \rightarrow \infty$) as $\gamma \rightarrow 0$. To observe this, let us assume to the contrary that there exists $\bb \in \tilde{H}^{1/2}(\Gamma_{\!\alpha} \cup {\Gamma}_{\!\circ})^3$ such that 
\beq \notag 
\begin{aligned}
\bPhi_L^\infty(\ba)(\hat\bxi) ~=~ -\Big(& \text{\emph{i}} k_p \,  \hat\bxi \exs \int_{\Gamma_{\!\alpha} \cup {\Gamma}_{\!\circ}} \, \Big\lbrace \lambda \exs (\bb \sip \bn) \,+\, 2\mu \exs (\bn \sip \hat\bxi) ( \bb \sip \hat\bxi)  \Big\rbrace  e^{-\text{\emph{i}}k_p \hat\bxi \cdot \by} \,\, \text{d}S_{\by} \\ 
 & \textcolor{black}{\oplus}~ \text{\emph{i}} k_s \,  \hat\bxi \exs \times \int_{\Gamma_{\!\alpha} \cup {\Gamma}_{\!\circ}} \Big\lbrace   \mu \exs(\bb \times \hat\bxi)(\bn \sip \hat\bxi) \,+\, \mu \exs  (\bn \times \hat\bxi) (\bb \sip \hat\bxi)  \Big\rbrace \, e^{-\textrm{\emph{i}} k_s \hat\bxi \cdot \by} \,\, \text{d}S_{\by} \Big), 
\end{aligned}
\eeq 
 associated with the layer potential
 \vspace{-1mm}
\beq\lb{Dlpb}
\bPhi_{\exs \Gamma_{\!\alpha} \cup {\Gamma}_{\!\circ}}(\bxi) ~=~ \int_{\Gamma_{\!\alpha} \cup {\Gamma}_{\!\circ}} \bb(\by) \cdot  \bfT(\bxi,\by)  \, \text{d}S_{\by}, \quad \bfT(\bxi,\by) ~=~ \bn(\by)\cdot\bSig(\bxi,\by), \quad \bxi \in \mathcal{B} \backslash \lbrace \Gamma_{\!\alpha} \cup {\Gamma}_{\!\circ} \rbrace.
\vspace{-1mm}
\eeq
On the other hand, owing to Definition of $\bPhi_L^\infty(\hat\bxi)$ in~\cref{Phi-inf}, the potential $\bPhi_{\exs \Gamma_{\!\alpha} \cup {\Gamma}_{\!\circ}}(\bxi)$ should coincide with 
 \vspace{-1mm}
\beq\lb{Pb}
\bPhi_{\! L}(\bxi) ~=~ \int_L \ba(\by) \cdot  \bfT(\bxi,\by)  \,\, \text{d}S_{\by},  \qquad \bxi \in \mathcal{B} \backslash L,
 \vspace{-1mm}
\eeq
over $ \bxi \in \mathcal{B} \backslash (L \cup \Gamma_{\!\alpha}  \cup {\Gamma}_{\!\circ})$. Now, let $\Gamma_{\!\alpha}  \cup {\Gamma}_{\!\circ} \not\ni\bxio\!\in L$ and let $\mathcal{B}_\epsilon$ be a small ball centered at $\bxio$ such that $\mathcal{B}_\epsilon \cap \lbrace \Gamma_{\!\alpha}  \cup {\Gamma}_{\!\circ} \rbrace = \emptyset $. In this case $\bPhi_{\exs \Gamma_{\!\alpha}  \cup {\Gamma}_{\!\circ}}\!$ is analytic in $\mathcal{B}_\epsilon$, while $\bPhi_{\! L}$ has a discontinuity across $\mathcal{B}_\epsilon\cap L$ -- which by contradiction completes the proof. 
\end{itemize}
\end{proof}

\vspace{-1mm}
\begin{theorem}[invariants of the solution sequence $(\bg_\alpha)_{\alpha \in {\mathbb{Z}}^*}\!$ for noise-free data] \lb{INV1}
Define 
\vspace{-1mm}
\beq \lb{Inv1}
\begin{aligned}
\Lambda_{\alpha}(\bg_\alpha, \bg_{\alpha+1}) &:= \big(\, \bg_{\alpha+1}\!-\bg_\alpha,\, F_{\alpha_\sharp}(\bg_{\alpha+1}\!-\bg_\alpha)\big), \\*[0.5mm]
\Upsilon_{\!\alpha}(\bg_\alpha, \bg_{\alpha+1}) &:= \big| \big(\, \bg_{\alpha+1},\, F_{{\alpha+1}_\sharp}\exs\bg_{\alpha+1}\big) - \big(\, \bg_{\alpha},\, F_{{\alpha}_\sharp}\exs\bg_{\alpha}\big) \big|, \quad \bg_\alpha, \bg_{\alpha+1}\in L^2(\OOd)^3,
\end{aligned}
\vspace{-1mm}
\eeq
where $(\bg_\alpha, \bg_{\alpha+1})(L;\gamma)$ are the constructed minimizers of $(\mathfrak{J}^{\gamma}_\alpha, \mathfrak{J}^{\gamma}_{\alpha+1})$ in~\eqref{GCf} according to~\eqref{mseq1}. Then, on defining $\tilde{\Gamma}_{\alpha+1}\! := {\Gamma}_{\!\alpha+1} \nxs\cap\nxs \lbrace \Gamma_{\!\alpha} \cup {\Gamma}_{\!\circ} \rbrace$, it follows that 
\vspace{1mm}
\begin{itemize}
\item{} If $L\subset\Gamma_{\!\alpha} \cup {\Gamma}_{\!\circ} \backslash \tilde{\Gamma}_{\alpha+1}$ then
\vspace{-3mm}
\beq\nonumber
\lim\limits_{\gamma \rightarrow 0}\Lambda_{\alpha}[\,\bg_\alpha, \bg_{\alpha+1}](L;\gamma) = 0 \,\, \text{and} \,\, \lim\limits_{\gamma \rightarrow 0}\Upsilon_{\!\alpha}[\,\bg_\alpha, \bg_{\alpha+1}](L;\gamma) = 0.
\vspace{-2mm}
\eeq
\item{} If $L\subset\tilde{\Gamma}_{\alpha+1}$ then
\vspace{-2mm}
\beq\nonumber
0 < \lim\limits_{\gamma \rightarrow 0}\Lambda_{\alpha}[\,\bg_\alpha, \bg_{\alpha+1}](L;\gamma) < \infty \,\, \text{and} \,\, 0 < \lim\limits_{\gamma \rightarrow 0}\Upsilon_{\!\alpha}[\,\bg_\alpha, \bg_{\alpha+1}](L;\gamma) < \infty. 
\vspace{-2mm}
\eeq
\item{} If $L\subset\hat{\Gamma}_{\!\alpha+1}\!$ then 
\vspace{-2mm}
\beq\nonumber
\lim\limits_{\gamma \rightarrow 0}\Lambda_{\alpha}[\,\bg_\alpha, \bg_{\alpha+1}](L;\gamma) = \infty \,\, \text{and} \,\,\lim\limits_{\gamma \rightarrow 0}\Upsilon_{\!\alpha}[\,\bg_\alpha, \bg_{\alpha+1}](L;\gamma) = \infty.
\vspace{-2mm}
\eeq
\end{itemize}
\end{theorem} 

\begin{proof}
If $L\subset\Gamma_{\!\circ} \cup \exs {\Gamma}_{\!\alpha} \cup \hat{\Gamma}_{\!\alpha+1}$, then on recalling that (a) $\mathcal{H}_{\alpha}^* \dbva = \bPhi_L^\infty$ (resp.~$\mathcal{H}_{\alpha+1}^* \dbvap = \bPhi_L^\infty$) from~\cref{TR1}, (b) continuity of $T_\alpha$ (resp.~$T_{\alpha+1}$) and its inverse $T_\alpha^{-1}$ (resp.~$T_{\alpha+1}^{-1}$) as per~\cref{Prelim}, and (c) the fact that $\bPhi_L^\infty \in Range(G_\alpha)$ (resp.~$\bPhi_L^\infty \in Range(G_{\alpha+1}))$ in~\cref{GLSM1}, it follows that 
\vspace*{-1mm}
\[
\lim\limits_{\gamma \rightarrow 0}\mathcal{H}_{\alpha}\bg_\alpha(L;\gamma) = T_{\alpha}^{-1} \dbva, (\text{resp.}~\lim\limits_{\gamma \rightarrow 0}\mathcal{H}_{\alpha+1}\bg_{\alpha+1}(L;\gamma) = T_{\alpha+1}^{-1} \dbvap), \,\,\, \bg_\alpha, \bg_{\alpha+1}\in L^2(\OOd)^3.
\vspace*{-1mm}
\] 
In this setting,
\vspace{1mm}
\begin{itemize}
\item{} If $L\subset\Gamma_{\!\circ} \cup {\Gamma}_{\!\alpha} \backslash \tilde{\Gamma}_{\alpha+1}$, then observe that (a) $\bK_{\!\alpha+1\!} = \bK_{\!\alpha}$ on $\Gamma_{\!\alpha}\backslash \tilde{\Gamma}_{\alpha+1}$ in the field equations~\eqref{GEi}, and (b) $\dbvap = \dbva$ on $\Gamma_{\!\alpha} \cup {\Gamma}_{\!\circ}$ and $\dbvap = \bzero$ on $\hat{\Gamma}_{\!\alpha+1}\!$ as per~\cref{TR1}. Accordingly, $\bv_{\alpha} = \bv_{\alpha+1}$ in $\mathcal{B} \backslash \lbrace \Gamma_{\!\alpha} \cup {\Gamma}_{\!\circ} \rbrace$ in light of the respective integral representations, $\forall \bxi \in \mathcal{B} \backslash \lbrace \Gamma_{\!\alpha} \cup {\Gamma}_{\!\circ} \rbrace$,
 \vspace{-1mm}
\beq\lb{Dlpb2}
\bv_{\alpha}(\bxi) = \int_{\Gamma_{\!\alpha} \cup {\Gamma}_{\!\circ}} \dbva(\by) \cdot  \bfT(\bxi,\by)  \, \text{d}S_{\by}, \,\, \bv_{\alpha+1}(\bxi) = \int_{\Gamma_{\!\alpha} \cup {\Gamma}_{\!\circ}} \dbvap(\by) \cdot  \bfT(\bxi,\by)  \, \text{d}S_{\by}, 
\vspace{-1mm}
\eeq
where $\bfT(\bxi,\by)$ is defined in~\eqref{Dlpb}. As a result, the contact laws in~\eqref{GEi}, governing $T_{\alpha}$ and $T_{\alpha}^{-1}$, read $T_{\alpha+1}^{-1}\dbvap =T_{\alpha}^{-1}\dbva$ and $T_{{\alpha+1}_\sharp}(\cdot)=T_{{\alpha}_\sharp}(\cdot)$ on $\Gamma_{\!\alpha} \cup {\Gamma}_{\!\circ}$. Thus, $\lim\limits_{\gamma \rightarrow 0}\mathcal{H}_{\alpha}\bg_\alpha(L;\gamma) = \lim\limits_{\gamma \rightarrow 0}\mathcal{H}_{\alpha}\bg_{\alpha+1}(L;\gamma) = T_{\alpha}^{-1} \dbva$ on $\Gamma_{\!\alpha} \cup {\Gamma}_{\!\circ}$. Invoking~\eqref{facts2}, one concludes 
\vspace*{-1mm}
\[
\begin{aligned}
\lim\limits_{\gamma \rightarrow 0}\Lambda_{\alpha}[\,\bg_\alpha, \bg_{\alpha+1}](L;\gamma) &= \lim\limits_{\gamma \rightarrow 0} \big(\, \mathcal{H}_{\alpha}\bg_{\alpha+1}\!-\mathcal{H}_{\alpha}\bg_\alpha,\, T_{\alpha_\sharp}(\mathcal{H}_{\alpha}\bg_{\alpha+1}\!-\mathcal{H}_{\alpha}\bg_\alpha)\big) = 0, \\*[0.5mm]
\lim\limits_{\gamma \rightarrow 0}\Upsilon_{\alpha}[\,\bg_\alpha, \bg_{\alpha+1}](L;\gamma) &=  \lim\limits_{\gamma \rightarrow 0} \big| \big(\, \mathcal{H}_{\alpha}\bg_{\alpha},\, (T_{{\alpha+1}_\sharp}-T_{{\alpha}_\sharp})\exs\mathcal{H}_{\alpha}\bg_{\alpha}\big) \big| = 0.
\end{aligned}
\vspace*{-1mm}
\]

\item{} If $L\subset\tilde{\Gamma}_{\alpha+1}$, then $\bK_{\!\alpha+1\!} \neq \bK_{\!\alpha}$ on $\tilde{\Gamma}_{\alpha+1\!}$ while $\dbvap = \dbva \in \tilde{H}^{1/2}(\Gamma_{\!\alpha} \cup {\Gamma}_{\!\circ})$, and $\dbvap = \bzero$ on $\hat{\Gamma}_{\!\alpha+1}\!$. Thus, $\lim\limits_{\gamma \rightarrow 0}\mathcal{H}_{\alpha}\bg_\alpha(L;\gamma) = T_{\alpha}^{-1} \dbva$ and $\lim\limits_{\gamma \rightarrow 0}\mathcal{H}_{\alpha+1}\bg_{\alpha+1}(L;\gamma) = \lim\limits_{\gamma \rightarrow 0}\mathcal{H}_{\alpha}\bg_{\alpha+1}(L;\gamma) = T_{\alpha+1}^{-1} \dbvap$. Note that here $T_{\alpha}^{-1}\dbva \in {H}^{-1/2}(\Gamma_{\!\alpha} \cup {\Gamma}_{\!\circ})$ and $T_{\alpha+1}^{-1} \dbvap \in {H}^{-1/2}(\Gamma_{\!\alpha+1} \cup {\Gamma}_{\!\circ})$ while $T_{\alpha}^{-1}\dbva \neq T_{\alpha+1}^{-1} \dbvap$ and $T_{{\alpha+1}_\sharp}(\cdot)\neq T_{{\alpha}_\sharp}(\cdot)$.  Therefore, 
\vspace*{-1mm}
\beq\nonumber
\begin{aligned}
\lim\limits_{\gamma \rightarrow 0}\Lambda_{\alpha}[\,\bg_\alpha, \bg_{\alpha+1}](L;\gamma) &=  \big((T_{\alpha+1}^{-1} -T_{\alpha}^{-1}) \dbva,\, T_{\alpha_\sharp}(T_{\alpha+1}^{-1} -T_{\alpha}^{-1}) \dbva\big) \\[0.15mm]
& \leqslant\, \norms{\!T_{\alpha_\sharp}\!}_{\tilde{H}^{\frac{1}{2}}(\Gamma_{\!\alpha} \cup {\Gamma}_{\!\circ})} \norms{\!(T_{\alpha+1}^{-1} -T_{\alpha}^{-1}) \dbva\!}^2_{{H}^{-\frac{1}{2}}(\Gamma_{\!\alpha} \cup {\Gamma}_{\!\circ})}  , 
\end{aligned}
\vspace*{-1mm}
\eeq

\vspace*{-1mm}
\[
\begin{aligned}
& \lim\limits_{\gamma \rightarrow 0}\Upsilon_{\alpha}[\,\bg_\alpha, \bg_{\alpha+1}](L;\gamma) =   \big| \big(\, T_{\alpha+1}^{-1} \dbva,\, T_{{\alpha+1}_\sharp}\exs T_{\alpha+1}^{-1} \dbva\big) - \big(\, T_{\alpha}^{-1} \dbva,\, T_{{\alpha}_\sharp}\exs T_{\alpha}^{-1} \dbva\big) \big| \\[0.15mm]
& \leqslant \, \big{(}\!\norms{\!T_{{\alpha+1}_\sharp}\!}_{\tilde{H}^{\frac{1}{2}}(\Gamma_{\!\alpha} \cup {\Gamma}_{\!\circ})} \norms{\!T_{\alpha+1}^{-1} \dbva\!}^2_{{H}^{-\frac{1}{2}}(\Gamma_{\!\alpha} \cup {\Gamma}_{\!\circ})} \!\!+\norms{\!T_{{\alpha}_\sharp}\!}_{\tilde{H}^{\frac{1}{2}}(\Gamma_{\!\alpha} \cup {\Gamma}_{\!\circ})} \norms{\!T_{\alpha}^{-1} \dbva\!}^2_{{H}^{-\frac{1}{2}}(\Gamma_{\!\alpha} \cup {\Gamma}_{\!\circ})} \!\big),
\end{aligned}
\vspace*{1mm}
\] 

\item{} If $L\subset\hat{\Gamma}_{\!\alpha+1}$, then $\lim\limits_{\gamma \rightarrow 0}\mathcal{H}_{\alpha+1}\bg_{\alpha+1}(L;\gamma) = T_{\alpha+1}^{-1} \dbvap \in {{H}^{-1/2}(\Gamma_{\!\alpha+1} \cup {\Gamma}_{\!\circ})}$ while $\lim\limits_{\gamma \rightarrow 0}\mathcal{H}_{\alpha}\bg_{\alpha}(L;\gamma) = \infty$. Invoking~\eqref{Inv1} and~\eqref{co-T}, one may observe
\vspace*{-2mm}
\beq\nonumber
\begin{aligned}
\Lambda_{\alpha}(\bg_\alpha, \bg_{\alpha+1}) = & \big(\,  \mathcal{H}_{\alpha}\bg_{\alpha+1}\!-\mathcal{H}_{\alpha}\bg_\alpha,\, T_{\alpha_\sharp}(\mathcal{H}_{\alpha}\bg_{\alpha+1}\!-\mathcal{H}_{\alpha}\bg_\alpha)\big) \\
& \geqslant \textrm{c} \nxs \norms{\!\mathcal{H}_{\alpha}\bg_{\alpha+1}\!-\mathcal{H}_{\alpha}\bg_\alpha\!}_{H^{-\frac{1}{2}}(\Gamma_{\!\circ} \cup \Gamma_\alpha)}^2,
\end{aligned}
\vspace*{-2mm}
\eeq
indicating that $\lim\limits_{\gamma \rightarrow 0}\Lambda_{\alpha}[\,\bg_\alpha, \bg_{\alpha+1}](L;\gamma) = \infty$. 
 
 A similar argument results in $\lim\limits_{\gamma \rightarrow 0}\Upsilon_{\!\alpha}[\,\bg_\alpha, \bg_{\alpha+1}](L;\gamma) = \infty$.
\end{itemize}
\vspace{-3mm}
\end{proof}

\begin{theorem}[invariants of the solution sequence for noisy data] \lb{INV2}
Define 
\vspace{-1mm}
\beq \lb{Inv2}
\begin{aligned}
\Lambda^\delta_{\alpha}(\bg_\alpha, \bg_{\alpha+1}) &:= \big(\, \bg_{\alpha+1}\!-\bg_\alpha,\, F^\delta_{\alpha_\sharp}(\bg_{\alpha+1}\!-\bg_\alpha)\big) + \exs\delta \norms{\bg_{\alpha+1}\!-\bg_\alpha\!}^2, \\*[0.5mm]
\Upsilon^\delta_{\!\alpha}(\bg_\alpha, \bg_{\alpha+1}) &:= \big| \Lambda^\delta_{\alpha+1}(\bg_{\alpha+1}, \bzero) -\Lambda^\delta_{\alpha}(\bg_\alpha, \bzero) \big|, \qquad \qquad \bg_\alpha, \bg_{\alpha+1}\in L^2(\OOd)^3,
\end{aligned}
\vspace{-1mm}
\eeq
where $(\bg_\alpha, \bg_{\alpha+1})(L;\delta,\gamma)$ are the constructed minimizers of $(\mathfrak{J}^{\delta,\gamma}_\alpha, \mathfrak{J}^{\delta,\gamma}_{\alpha+1})$ in~\eqref{GCfn} according to~\eqref{gL}. Then,
\vspace{1mm}
\begin{itemize}
\item{} If $L\subset\Gamma_{\!\alpha} \cup {\Gamma}_{\!\circ} \backslash \tilde{\Gamma}_{\alpha+1}$ then
\vspace*{-2mm}
\beq\nonumber
\lim\limits_{\gamma \rightarrow 0}\liminf\limits_{\delta \rightarrow 0}\Lambda_{\alpha}^\delta[\,\bg_\alpha, \bg_{\alpha+1}](L;\delta,\gamma) = 0 \,\, \text{and} \,\,\lim\limits_{\gamma \rightarrow 0}\liminf\limits_{\delta \rightarrow 0}\Upsilon^\delta_{\!\alpha}[\,\bg_\alpha, \bg_{\alpha+1}](L;\delta,\gamma) = 0.
\vspace*{-2mm}
\eeq
\item{} If $L\subset\tilde{\Gamma}_{\alpha+1}$ then
\vspace*{-2mm}
\beq\nonumber
0\nxs<\nxs\lim\limits_{\gamma \rightarrow 0}\liminf\limits_{\delta \rightarrow 0}\Lambda_{\alpha}^\delta[\,\bg_\alpha, \bg_{\alpha+1}](L;\delta,\gamma) \!<\! \infty \, \text{and} \, 0\nxs<\nxs\lim\limits_{\gamma \rightarrow 0}\liminf\limits_{\delta \rightarrow 0}\Upsilon^\delta_{\!\alpha}[\,\bg_\alpha, \bg_{\alpha+1}](L;\delta,\gamma) \!< \! \infty.
\vspace*{-2mm}
\eeq
\item{} If $L\subset\hat{\Gamma}_{\!\alpha+1}\!$ then 
\vspace*{-2mm}
\beq\nonumber
\lim\limits_{\gamma \rightarrow 0}\liminf\limits_{\delta \rightarrow 0}\Lambda_{\alpha}^\delta[\,\bg_\alpha, \bg_{\alpha+1}](L;\delta,\gamma) = \infty\, \text{and} \,\lim\limits_{\gamma \rightarrow 0}\liminf\limits_{\delta \rightarrow 0}\Upsilon^\delta_{\!\alpha}[\,\bg_\alpha, \bg_{\alpha+1}](L;\delta,\gamma) = \infty.
\vspace*{-2mm}
\eeq
\end{itemize}
\end{theorem} 
\vspace{-1mm}
\begin{proof}
If $L \subset \Gamma_{\!\alpha} \cup {\Gamma}_{\!\circ}$, then~\cref{GLSM2} and~\cref{TR1} indicate that there exists a sequence $\delta(\gamma)$ such that 
\vspace*{-2mm}
\beq\nonumber
\limsup\limits_{\gamma \rightarrow 0} \norms{\nxs\mathcal{H}_{\alpha}\bg_{\alpha}^{\delta(\gamma), \gamma}\!-T_{\alpha}^{-1}\dbva\nxs}^2 \,=\, 0, \quad \limsup\limits_{\gamma \rightarrow 0} \norms{\nxs\mathcal{H}_{\alpha+1}\exs\bg_{\alpha+1}^{\delta(\gamma), \gamma}\!-T_{\alpha+1}^{-1}\dbvap\nxs}^2 \,=\, 0.
\vspace*{-2mm}
\eeq
From~\cref{INV1}, one may then conclude that for $L\subset\Gamma_{\!\alpha} \cup {\Gamma}_{\!\circ} \backslash \tilde{\Gamma}_{\alpha+1}$, 
\vspace*{-2mm}
\beq\nonumber
\lim\limits_{\gamma \rightarrow 0}\Lambda_{\alpha}[\,\bg_\alpha, \bg_{\alpha+1}](L;\delta(\gamma),\gamma) = 0, \quad \lim\limits_{\gamma \rightarrow 0}\Upsilon_{\!\alpha}[\,\bg_\alpha, \bg_{\alpha+1}](L;\delta(\gamma),\gamma) = 0, 
\vspace*{-2mm}
\eeq  
and for $L\subset\tilde{\Gamma}_{\alpha+1}$,   
\vspace*{-2mm}
\beq\nonumber
0 < \lim\limits_{\gamma \rightarrow 0}\Lambda_{\alpha}[\,\bg_\alpha, \bg_{\alpha+1}](L;\delta(\gamma),\gamma) < \infty , \quad 0 < \lim\limits_{\gamma \rightarrow 0}\Upsilon_{\!\alpha}[\,\bg_\alpha, \bg_{\alpha+1}](L;\delta(\gamma),\gamma) < \infty. 
\vspace*{-2mm}
\eeq
One may also observe that
\vspace{-1mm}
\beq\nonumber
\begin{aligned}
&\big| \big[ \Lambda^\delta_{\alpha} -\Lambda_{\alpha} \big] (\bg_\alpha, \bg_{\alpha+1}) \big| \,\leqslant\, 2 \delta \big{(}\!\norms{\!\bg_\alpha\!}^2 \nxs +\nxs \norms{\!\bg_{\alpha+1}\!}^2\!\big{)}, \\[0.1mm]
&\big| \big[ \Upsilon^\delta_{\alpha} -\Upsilon_{\alpha} \big] (\bg_\alpha, \bg_{\alpha+1}) \big| \,\leqslant\, 2 \delta \big{(}\!\norms{\!\bg_\alpha\!}^2 \nxs +\nxs \norms{\!\bg_{\alpha+1}\!}^2\!\big{)}. 
\end{aligned}
\vspace{-1mm}
\eeq

Invoking~\eqref{lmt}, one then deduces that
\vspace{-2mm}
\beq\nonumber
\begin{aligned}
& \limsup\limits_{\gamma \rightarrow 0} \exs \big| \big( \Lambda^{\delta(\gamma)}_{\alpha} -\Lambda_{\alpha} \big) [\bg_\alpha, \bg_{\alpha+1}](L;\delta(\gamma),\gamma) \big| ~\!\!=\!~ 0, \\[-0.1mm]
& \limsup\limits_{\gamma \rightarrow 0} \exs \big| \big( \Upsilon^{\delta(\gamma)}_{\alpha} -\Upsilon_{\alpha} \big) [\bg_\alpha, \bg_{\alpha+1}](L;\delta(\gamma),\gamma) \big| ~\!\!=\!~ 0,
\end{aligned}
\vspace{-3mm}
\eeq
which completes the proof for the first two parts of the theorem. If $L\subset\hat{\Gamma}_{\!\alpha+1}\!$, then based on \eqref{co-T}, \eqref{Ns-op}, \eqref{statG}, and \eqref{statG2}, one may find that 
 \vspace*{-1.5mm}
\[
\begin{aligned}
&\Lambda_{\alpha}^{\delta}(\bg_\alpha, \bg_{\alpha+1}) \geqslant \textrm{c}_{\circ} \nxs \norms{\!\mathcal{H}_{\alpha}\bg_{\alpha+1}\!-\mathcal{H}_{\alpha}\bg_\alpha\!}_{H^{-\frac{1}{2}}(\Gamma_{\!\circ} \cup \Gamma_\alpha)}^2, \\*[0.1mm]
&\Upsilon_{\alpha}^{\delta}(\bg_\alpha, \bg_{\alpha+1}) \geqslant \big{|} \textrm{c}_{1}\norms{\! \mathcal{H}_{\alpha}\bg_{\alpha}\!}_{H^{-\frac{1}{2}}(\Gamma_{\!\circ} \cup \Gamma_{\alpha})}^2\!\!- \, \textrm{c}_{2}\norms{\!\mathcal{H}_{\alpha+1}\bg_{\alpha+1}\!}_{H^{-\frac{1}{2}}(\Gamma_{\!\circ} \cup \Gamma_{\alpha+1})}^2 \big{|},
\end{aligned}
\vspace*{-1.5mm}
\] 
where $\textrm{c}_\circ, \textrm{c}_1, \textrm{c}_2 > 0$ are constants independent of $\mathcal{H}_{\alpha}\bg_\alpha$. The theorem's statement will then follow in light of  
\vspace{-1mm}
\[
\limsup\limits_{\gamma \rightarrow 0}\limsup\limits_{\delta \rightarrow 0} \norms{\!\mathcal{H}_{\alpha+1}\exs\bg_{\alpha+1}\!}_{H^{-\frac{1}{2}}(\Gamma_{\!\circ} \cup \Gamma_{\alpha+1})}^2 < \infty, \quad \liminf\limits_{\gamma \rightarrow 0}\liminf\limits_{\delta \rightarrow 0} \norms{\!\mathcal{H}_{\alpha}\exs\bg_{\alpha}\!}_{H^{-\frac{1}{2}}(\Gamma_{\!\circ} \cup \Gamma_\alpha)}^2 = \infty.
\vspace{-1mm}
\]
\end{proof}

Based on Theorems 4.5-7, the differential evolution indicators are introduced in the sequel.

 {\bfseries\slshape Differential evolution indicators for noise-free data.}~Let us introduce the evolution indicator functionals $I_\alpha^\mathcal{D}: L^2(\Omega^3)\times L^2(\Omega^3) \rightarrow \mathbb{R}$ and $\hat{I}_\alpha^\mathcal{D}: L^2(\Omega^3)\times L^2(\Omega^3) \rightarrow \mathbb{R}$ such that 
\vspace*{-1mm}
\beq\lb{EIF}
\begin{aligned}
& I_\alpha^\mathcal{D}(\bg_\alpha, \bg_{\alpha+1}) := \frac{1}{\sqrt{\Lambda_{\alpha+1}(\bzero, \bg_{\alpha+1}) \big{[}1+ \Lambda_{\alpha+1}(\bzero, \bg_{\alpha+1})\mathcal{D}_{\alpha}^{-1}(\bg_\alpha, \bg_{\alpha+1})\big{]}}},  \\*[0.25mm]
&\hat{I}_\alpha^\mathcal{D}(\bg_\alpha, \bg_{\alpha+1}) := \frac{1}{\sqrt{\Lambda_{\alpha}(\bg_{\alpha}, \bzero) + \Lambda_{\alpha+1}(\bzero, \bg_{\alpha+1}) \big{[}1+ \Lambda_{\alpha}(\bg_{\alpha}, \bzero)\mathcal{D}_{\alpha}^{-1}(\bg_\alpha, \bg_{\alpha+1})\big{]}}}, 
\end{aligned}
\vspace*{-1mm}
\eeq
where $\mathcal{D}_{\alpha} \in \lbrace \Lambda_{\alpha}, \Upsilon_{\alpha} \rbrace$, and $(\bg_\alpha, \bg_{\alpha+1})(L;\gamma) \in L^2(\Omega^3)\times L^2(\Omega^3)$ are the constructed minimizers of $(\mathfrak{J}^{\gamma}_\alpha, \mathfrak{J}^{\gamma}_{\alpha+1})$ in~\eqref{GCf} according to~\eqref{mseq1}.
Then, it follows that
\vspace*{1mm}
\begin{itemize}
\item{} $L \subset \tilde{\Gamma}_{\alpha+1} \cup \hat{\Gamma}_{\alpha+1} \iff \lim\limits_{\gamma \rightarrow 0} I_\alpha^\mathcal{D}(\bg_\alpha, \bg_{\alpha+1})(L;\gamma) > 0$. 
\item{} $L \subset  \hat{\Gamma}_{\alpha+1} \iff \lim\limits_{\gamma \rightarrow 0} \hat{I}_\alpha^\mathcal{D}(\bg_\alpha, \bg_{\alpha+1})(L;\gamma) > 0$. 
\end{itemize}

This may be observed (a) by invoking~\cref{GLSM1} which reads $L \subset \Gamma_{\!\circ} \cup \Gamma_{\alpha+1}$ (\emph{resp.} $L \subset \Gamma_{\!\circ} \cup \Gamma_{\alpha}$) if and only if $\lim\limits_{\gamma \rightarrow 0} \Lambda_{\alpha+1}(\bzero, \bg_{\alpha+1})(L;\gamma) < \infty$ (\emph{resp.} $\lim\limits_{\gamma \rightarrow 0} \Lambda_{\alpha}(\bg_{\alpha}, \bzero)(L;\gamma) < \infty$), implying that $\lim\limits_{\gamma \rightarrow 0} I_\alpha^\mathcal{D}(\bzero, \bg_{\alpha+1})(L;\gamma)~=~0$ (\emph{resp.} $\lim\limits_{\gamma \rightarrow 0} \hat{I}_\alpha^\mathcal{D}(\bg_{\alpha}, \bzero)(L;\gamma) = 0$) when $L \subset \mathcal{B} \backslash \lbrace \Gamma_{\!\circ} \cup \Gamma_{\alpha+1} \rbrace$ (\emph{resp.} $L \subset \mathcal{B} \backslash \lbrace \Gamma_{\!\circ} \cup \Gamma_{\alpha} \rbrace$), (b) in view of the first statement of~\cref{INV1} which ensures that for $L \subset \Gamma_{\!\circ} \cup \Gamma_{\alpha} \backslash \tilde{\Gamma}_{\alpha+1}$, $\lim\limits_{\gamma \rightarrow 0} I_\alpha^\mathcal{D}(\bg_\alpha, \bg_{\alpha+1})(L;\gamma) = 0$ and $\lim\limits_{\gamma \rightarrow 0} \hat{I}_\alpha^\mathcal{D}(\bg_\alpha, \bg_{\alpha+1})(L;\gamma) = 0$, and (c) by recalling the second and third statements of~\cref{INV1}. In other words, $\hat{I}_\alpha^\mathcal{D}$ illuminates the support of \emph{geometric} evolution between $[t_{\alpha}, t_{\alpha+1}]$ by achieving its highest values at the loci of newly born interfaces $\hat{\Gamma}_{\alpha+1}$. However, ${I}_\alpha^\mathcal{D}$ more holistically reconstructs the support of the micromechanical evolution which includes the new interstitial spaces $\hat{\Gamma}_{\alpha+1}$ as well as the pre-existing interfaces $\tilde{\Gamma}_{\alpha+1}$ whose elastic properties have changed between $[t_{\alpha}, t_{\alpha+1}]$ e.g.~due to chemical reaction or micro-slip.      

 {\bfseries\slshape Differential evolution indicators for noisy data.}~Consider the indicator functionals $I_{\alpha,\delta}^\mathcal{D^\delta}: L^2(\Omega^3)\times L^2(\Omega^3) \rightarrow \mathbb{R}$ and $\hat{I}_{\alpha,\delta}^\mathcal{D^\delta}: L^2(\Omega^3)\times L^2(\Omega^3) \rightarrow \mathbb{R}$ such that 
\vspace*{-1mm}
\beq\lb{EIFn}
\begin{aligned}
& I_{\alpha,\delta}^\mathcal{D^\delta}(\bg_\alpha, \bg_{\alpha+1}) := \frac{1}{\sqrt{\Lambda^\delta_{\alpha+1}(\bzero, \bg_{\alpha+1}) \big{[}1+ \Lambda^\delta_{\alpha+1}(\bzero, \bg_{\alpha+1})\mathcal{D}_{\alpha,\delta}^{-1}(\bg_\alpha, \bg_{\alpha+1})\big{]}}},  \\
&\hat{I}_{\alpha,\delta}^\mathcal{D^\delta}(\bg_\alpha, \bg_{\alpha+1}) := \frac{1}{\sqrt{\Lambda^\delta_{\alpha}(\bg_{\alpha}, \bzero) + \Lambda^\delta_{\alpha+1}(\bzero, \bg_{\alpha+1}) \big{[}1+ \Lambda^\delta_{\alpha}(\bg_{\alpha}, \bzero)\mathcal{D}_{\alpha,\delta}^{-1}(\bg_\alpha, \bg_{\alpha+1})\big{]}}},
\end{aligned}
\vspace*{-1mm}
\eeq
where $\mathcal{D}_{\alpha,\delta} \in \lbrace \Lambda_{\alpha}^\delta, \Upsilon_{\alpha}^\delta \rbrace$, and $(\bg_\alpha, \bg_{\alpha+1})(L;\delta,\gamma) = (\bg_\alpha^{\delta,\gamma}, \bg_{\alpha+1}^{\delta,\gamma})  \in L^2(\Omega^3)\times L^2(\Omega^3)$ are the constructed minimizers of $(\mathfrak{J}^{\delta,\gamma}_\alpha, \mathfrak{J}^{\delta,\gamma}_{\alpha+1})$ in~\eqref{GCfn} according to~\eqref{gL}. 
Then, it follows that
\vspace*{1mm}
\begin{itemize}
\item{} $L \subset \tilde{\Gamma}_{\alpha+1} \cup \hat{\Gamma}_{\alpha+1} \iff \lim\limits_{\gamma \rightarrow 0} \exs \liminf\limits_{\delta(\gamma) \rightarrow 0} \exs I_{\alpha,\delta}^\mathcal{D^\delta} \exs (\bg^{\delta,\gamma}_\alpha, \bg^{\delta,\gamma}_{\alpha+1}) \,>\, 0$. 
\item{} $L \subset  \hat{\Gamma}_{\alpha+1} \iff \lim\limits_{\gamma \rightarrow 0} \exs \liminf\limits_{\delta(\gamma) \rightarrow 0} \exs \hat{I}_{\alpha,\delta}^\mathcal{D^\delta} \exs (\bg^{\delta,\gamma}_\alpha, \bg^{\delta,\gamma}_{\alpha+1}) \,>\, 0$. 
\end{itemize}

This may be established on the basis of (a)~\cref{GLSM2} which reads $L \subset \Gamma_{\!\circ} \cup \Gamma_{\alpha+1}$ (\emph{resp.} $L \subset \Gamma_{\!\circ} \cup \Gamma_{\alpha}$) if and only if $\lim\limits_{\gamma \rightarrow 0} \exs \liminf\limits_{\delta(\gamma) \rightarrow 0} \Lambda^\delta_{\alpha+1}(\bzero, \bg^{\delta,\gamma}_{\alpha+1}) < \infty$ (\emph{resp.} $\lim\limits_{\gamma \rightarrow 0} \exs \liminf\limits_{\delta(\gamma) \rightarrow 0} \Lambda^\delta_{\alpha}(\bg^{\delta,\gamma}_{\alpha}, \bzero) < \infty$), implying that $\lim\limits_{\gamma \rightarrow 0} \exs \liminf\limits_{\delta(\gamma) \rightarrow 0} I_{\alpha,\delta}^\mathcal{D^\delta}(\bzero, \bg^{\delta,\gamma}_{\alpha+1})~=~0$ (\emph{resp.} $\lim\limits_{\gamma \rightarrow 0} \exs \liminf\limits_{\delta(\gamma) \rightarrow 0} \hat{I}_{\alpha,\delta}^\mathcal{D^\delta}(\bg^{\delta,\gamma}_{\alpha}, \bzero) = 0$) when $L \subset \mathcal{B} \backslash \lbrace \Gamma_{\!\circ} \cup \Gamma_{\alpha+1} \rbrace$ (\emph{resp.} $L \subset \mathcal{B} \backslash \lbrace \Gamma_{\!\circ} \cup \Gamma_{\alpha} \rbrace$), (b) first statement of~\cref{INV2} which ensures that for $L \subset \Gamma_{\!\circ} \cup \Gamma_{\alpha} \backslash \tilde{\Gamma}_{\alpha+1}$, $\lim\limits_{\gamma \rightarrow 0} \exs \liminf\limits_{\delta(\gamma) \rightarrow 0} I_{\alpha,\delta}^\mathcal{D^\delta}(\bg^{\delta,\gamma}_{\alpha}, \bg^{\delta,\gamma}_{\alpha+1})~=~0$ and $\lim\limits_{\gamma \rightarrow 0} \exs \liminf\limits_{\delta(\gamma) \rightarrow 0} \hat{I}_{\alpha,\delta}^\mathcal{D^\delta}(\bg^{\delta,\gamma}_{\alpha}, \bg^{\delta,\gamma}_{\alpha+1}) = 0$, and (c) second and third statements of~\cref{INV2}. In other words, $\hat{I}_{\alpha,\delta}^\mathcal{D^\delta}$ illuminates the support of \emph{geometric} evolution between $[t_{\alpha}, t_{\alpha+1}]$ by achieving its highest values at the loci of newly born interfaces $\hat{\Gamma}_{\alpha+1}$. However, $I_{\alpha,\delta}^\mathcal{D^\delta} $ reconstructs the support of the evolution more comprehensively including the new interfacial spaces $\hat{\Gamma}_{\alpha+1}$ as well as the pre-existing interfaces $\tilde{\Gamma}_{\alpha+1}$ whose elastic properties have changed between $[t_{\alpha}, t_{\alpha+1}]$.

\section{Implementation and results}\lb{IR}   

To illustrate the theoretical developments, this section examines the performance of differential evolution indicators~\eqref{EIF} and~\eqref{EIFn} through a set of numerical experiments and compares the results to those obtained by the generalized linear sampling method~\cite{Fatemeh2017}. In what follows the synthetic sensory data, namely the scattered fields $\bv_\alpha$ at sensing steps $t _{\alpha}= \lbrace t_{\circ}, t_1, t_2, ... \rbrace$, are simulated by a computational platform based on the elastodynamic boundary integral equations, see~\cite{Fatemeh2015,Fatemeh2016,Bon1999} for details of the computational method.

\vspace*{-2mm}
\subsection{Testing configuration}

Two test setups are considered as illustrated in~\cref{SetNum1} and~\cref{SetNum2} where an elastic plate of dimensions $3$ $\!\times\!$ $3$ $\!\times\!$ $0.02$ is endowed with (I) a randomly cracked damage zone, and (II) a pore zone. The shear modulus, mass density, and Poisson's ratio of the plate are taken as $\mu = 1$, $\rho = 1$ and $\nu = 0.25$, whereby the shear and compressional wave speeds read $c_s = 1$ and $c_p = 1.73$. In \emph{Setup I}, shown in~\cref{SetNum1}, the damage zone is comprised of randomly distributed cracks $\Gamma_{1}-\Gamma_{24}$ evolving hidden within the thickness of the specimen in seven time steps $t_1-t_7$. A detailed description of scatterers including the center $(x_c, y_c)$, length $\ell$, and orientation $\phi$ (with respect to $x$ axis) of each crack $\Gamma_{\kappa}$, $\kappa = \lbrace 1, 2, ..., 24 \rbrace$ is provided in~\cref{Num1}. All fractures in this configuration are traction-free i.e., the interfacial stiffness~$\bK(\bxi) = \bzero$ on $\bxi \in \!\!{\textstyle \bigcup\limits_{\kappa =1}^{24} \!\! \Gamma_{\!\kappa}}$. In \emph{Setup II}, depicted in~\cref{SetNum2}, a bubble zone is growing within the plate thickness, comprised of randomly distributed pores $\Pi_{1}-\Pi_{21}$ developing in seven time steps $t_1-t_7$. A detailed description of the specimen including the center $(x_c, y_c)$ and radius $r$ of each pore $\Pi_{\kappa}$, $\kappa = \lbrace 1, 2, ..., 21 \rbrace$ is provided in~\cref{Num2}.

\begin{figure}[!tp]
\center\includegraphics[width=0.85\linewidth]{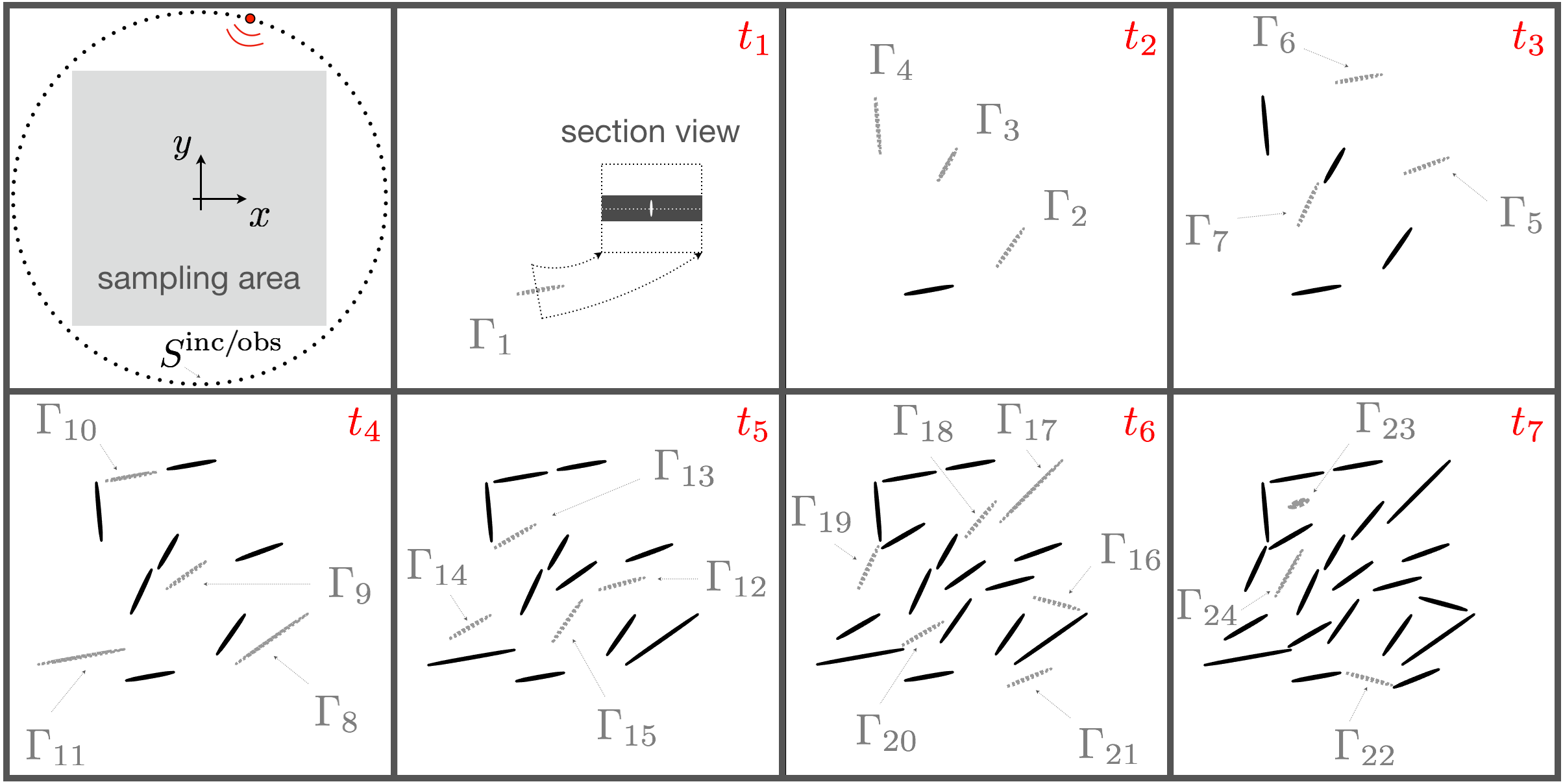} \vspace*{-3mm} 
\caption{Sensing configuration of synthetic sequential experiments on an elastic plate (top left), featuring a damage zone comprised of randomly distributed cracks $\Gamma_{1}-\Gamma_{24}$ evolving in seven time steps ($t_1-t_7$) within the thickness of the specimen according to the sectional view shown at $t_1$. } \lb{SetNum1}
\vspace*{-1mm}
\end{figure} 

\begin{table}[!h]
\vspace*{-1mm}
\begin{center}
\caption{\small Damage zone configuration illustrated in~\cref{SetNum1}:~center $(x_c, y_c)$, length $\ell$, and orientation $\phi$ (with respect to $x$ axis) of cracks $\Gamma_{\kappa}$, $\kappa = \lbrace 1, 2, ..., 24 \rbrace$.} \vspace*{-5mm}
\label{Num1}
\resizebox{\textwidth}{!}{%
\begin{tabular}{|c|c|c|c|c|c|c|c|c|c|c|c|c|} \hline
\!\!$\kappa$\!\! & 1 & 2 & 3 & 4&5 & 6 & 7 & 8 & 9 & 10 & 11 & 12 \\ \hline\hline 
\!\!$x_{\text{c}}(\Gamma_\kappa)$\!\!    & \!\!$-0.33$\!\!  & \!\!$0.21$\!\!  & \!\!$-0.21$\!\! & \!\!$-0.68$\!\! & \!\!$0.4$\!\! & \!\!$-0.05$\!\! & \!\!$-0.39$\!\!  & \!\!$0.49$\!\! & \!\!$-0.09$\!\! & \!\!$-0.46$\!\! & \!\!$-0.8$\!\! & \!\!$0.21$\!\! \\ 
 \hline
 \!\!$y_{\text{c}}(\Gamma_\kappa)$\!\!    & \!\!$-0.62$\!\!  & \!\!$-0.34$\!\!  & \!\!$0.22$\!\! & \!\!$0.49$\!\! & \!\!$0.21$\!\! & \!\!$0.8$\!\! & \!\!$-0.05$\!\!  & \!\!$-0.37$\!\! & \!\!$0.06$\!\! & \!\!$0.72$\!\! & \!\!$-0.5$\!\! & \!\!$0$\!\! \\ 
 \hline
\!\!$\ell \exs (\Gamma_\kappa)$\!\!   & \!\!$1/3$\!\!  & \!\!$1/3$\!\! & \!\!$1/4$\!\! & \!\!$2/5$\!\! & \!\!$1/3$\!\! & \!\!$1/3$\!\! & \!\!$1/3$\!\!  &  \!\!$3/5$\!\! & \!\!$1/3$\!\! & \!\!$2/5$\!\! & \!\!$3/5$\!\! & \!\!$1/3$\!\! \\ 
 \hline
 \!\!$\phi \exs (\Gamma_\kappa)$\!\!   & \!\!$\pi/18$\!\!  & \!\!$11\pi/36$\!\! & \!\!$\pi/3$\!\! & \!\!$19\pi/36$\!\! & \!\!$\pi/9$\!\!  & \!\!$\pi/18$\!\!  & \!\!$13\pi/36$\!\!  &  \!\!$7\pi/36$\!\! & \!\!$7\pi/36$\!\! & \!\!$\pi/18$\!\! & \!\!$\pi/18$\!\! & \!\!$\pi/12$\!\! \\ 
 \hline\hline
\!\!$\kappa$\!\! & 13 & 14 & 15 & 16&17 & 18 & 19 & 20 & 21 & 22 & 23 & 24 \\ \hline\hline 
\!\!$x_{\text{c}}(\Gamma_\kappa)$\!\!    & \!\!$-0.5$\!\!  & \!\!$-0.8$\!\!  & \!\!$-0.15$\!\! & \!\!$0.52$\!\! & \!\!$0.36$\!\! & \!\!$0.01$\!\! & \!\!$-0.74$\!\!  & \!\!$-0.38$\!\! & \!\!$0.34$\!\! & \!\!$0.02$\!\! & \!\!$-0.45$\!\! & \!\!$-0.51$\!\! \\ 
 \hline
 \!\!$y_{\text{c}}(\Gamma_\kappa)$\!\!    & \!\!$0.32$\!\!  & \!\!$-0.29$\!\!  & \!\!$-0.25$\!\! & \!\!$-0.13$\!\! & \!\!$0.62$\!\! & \!\! $0.43$\!\! & \!\!$0.1$\!\!  & \!\!$-0.34$\!\! & \!\!$-0.63$\!\! & \!\!$-0.64$\!\! & \!\!$0.55$\!\! & \!\!$0.08$\!\! \\ 
 \hline
\!\!$\ell \exs (\Gamma_\kappa)$\!\!   & \!\!$1/3$\!\!  & \!\!$1/3$\!\! & \!\!$7/20$\!\! & \!\!$1/3$\!\! & \!\!$3/5$\!\! & \!\!$1/3$\!\! & \!\!$1/3$\!\!  &  \!\!$1/3$\!\! & \!\!$1/3$\!\! & \!\!$1/3$\!\! & \!\!$1/7$\!\! & \!\!$1/3$\!\! \\ 
 \hline
 \!\!$\phi \exs (\Gamma_\kappa)$\!\!   & \!\!$\pi/6$\!\!  & \!\!$\pi/6$\!\! & \!\!$11\pi/36$\!\! & \!\!$-\pi/12$\!\! & \!\!$\pi/4$\!\!  & \!\!$5\pi/18$\!\!  & \!\!$13\pi/36$\!\!  &  \!\!$\pi/6$\!\! & \!\!$11\pi/90$\!\! & \!\!$-\pi/12$\!\! & \!\!$\pi/9$\!\! & \!\!$\pi/3$\!\! \\ 
 \hline
\end{tabular}}
\end{center}
\vspace*{-1mm}
\end{table}

\begin{figure}[!tp]
\center\includegraphics[width=0.85\linewidth]{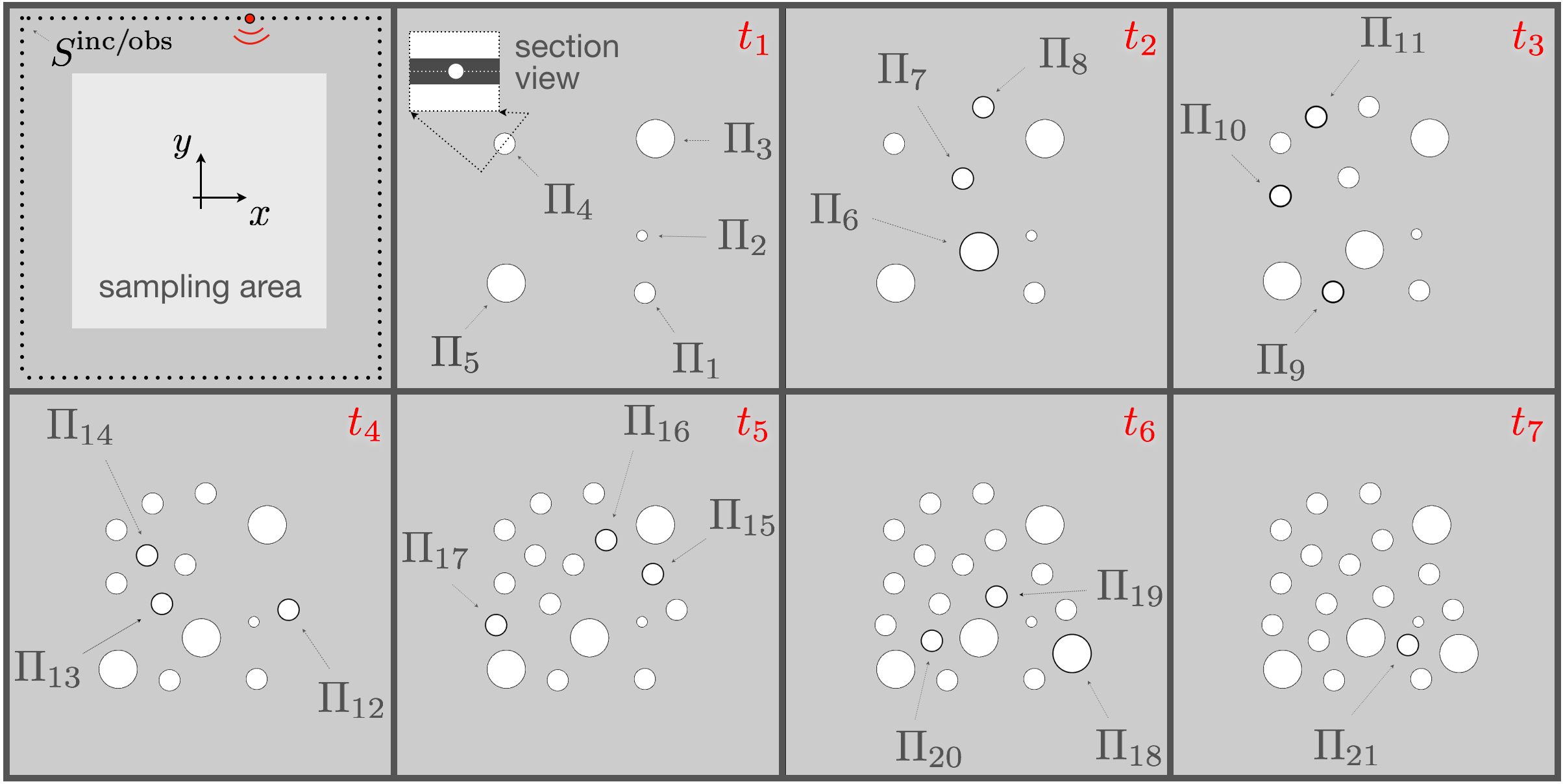} \vspace*{-3mm} 
\caption{Sensing configuration of synthetic sequential experiments on an elastic plate (top left), featuring a bubble zone comprised of randomly distributed pores $\Pi_{1}-\Pi_{21}$ evolving in seven time steps ($t_1-t_7$) within the thickness of the specimen according to the sectional view shown at $t_1$.} \lb{SetNum2}
\end{figure} 

\begin{table}[ht]
\vspace*{-1mm}
\begin{center}
\caption{\small Bubble zone configuration illustrated in~\cref{SetNum2}:~center $(x_c, y_c)$ and radius $r$ of bubbles $\Pi_{\kappa}$, $\kappa = \lbrace 1, 2, ..., 21 \rbrace$.}\vspace*{-3mm}
\label{Num2}
\begin{tabular}{|c|c|c|c|c|c|c|c|c|c|c|c|c|} \hline
\!\!$\kappa$\!\! & 1 & 2 & 3 & 4&5 & 6 & 7 & 8 & 9 & 10 & 11  \\ \hline\hline 
\!\!$x_{\text{c}}(\Pi_\kappa)$\!\!    & \!\!$0.42$\!\!  & \!\!$0.4$\!\!  & \!\!$0.51$\!\! & \!\!$-0.66$\!\! & \!\!$-0.65$\!\! & \!\!$0$\!\! & \!\!$-0.13$\!\!  & \!\!$0.03$\!\! & \!\!$-0.25$\!\! & \!\!$-0.66$\!\! & \!\!$-0.38$\!\!  \\ 
 \hline
 \!\!$y_{\text{c}}(\Pi_\kappa)$\!\!    & \!\!$-0.72$\!\!  & \!\!$-0.27$\!\!  & \!\!$0.48$\!\! & \!\!$0.44$\!\! & \!\!$-0.64$\!\! & \!\!$-0.4$\!\! & \!\!$0.17$\!\!  & \!\!$0.72$\!\! & \!\!$-0.73$\!\! & \!\!$0.02$\!\! & \!\!$0.64$\!\!  \\ 
 \hline
\!\!$r \exs (\Pi_\kappa)$\!\!   & \!\!$\exs 0.08$\!\!  & \!\!$\exs 0.04$\!\! & \!\!$\exs 0.15$\!\! & \!\!$\exs 0.08$\!\! & \!\!$\exs 0.15$\!\! & \!\!$\exs 0.15$\!\! & \!\!$\exs 0.08$\!\!  &  \!\!$0.08$\!\! & \!\!$\exs 0.08$\!\! & \!\!$\exs 0.08$\!\! & \!\!$\exs 0.08$\!\!  \\ 
 \hline\hline
\!\!$\kappa$\!\! & 12 & 13 & 14 & 15 &16 & 17 & 18 & 19 & 20 & 21 &   \\ \hline\hline 
\!\!$x_{\text{c}}(\Pi_\kappa)$\!\!    & \!\!$0.67$\!\!  & \!\!$-0.31$\!\!  & \!\!$-0.42$\!\! & \!\!$0.49$\!\! & \!\!$0.12$\!\! & \!\!$-0.73$\!\! & \!\!$0.72$\!\!  & \!\!$0.13$\!\! & \!\!$-0.37$\!\! & \!\!$0.32$\!\! & \!\!$ $\!\!  \\ 
 \hline
 \!\!$y_{\text{c}}(\Pi_\kappa)$\!\!    & \!\!$-0.18$\!\!  & \!\!$-0.13$\!\!  & \!\!$0.24$\!\! & \!\!$0.09$\!\! & \!\!$0.36$\!\! & \!\! $-0.3$\!\! & \!\!$-0.52$\!\!  & \!\!$-0.08$\!\! & \!\!$-0.42$\!\! & \!\!$-0.45$\!\! & \!\!$ $\!\!   \\ 
 \hline
\!\!$r \exs (\Pi_\kappa)$\!\!   & \!\!$\exs 0.08$\!\!  & \!\!$\exs 0.08$\!\! & \!\!$\exs 0.08$\!\! & \!\!$\exs 0.08$\!\! & \!\!$\exs 0.08$\!\! & \!\!$\exs 0.08$\!\! & \!\!$\exs 0.15$\!\!  &  \!\!$\exs 0.08$\!\! & \!\!$\exs 0.08$\!\! & \!\!$\exs 0.08$\!\! & \!\!$ $\!\!   \\ 
 \hline
\end{tabular}
\end{center}
\vspace*{-6mm}
\end{table}

\vspace*{-2mm}
\subsection{Forward scattering simulations}
Numerical experiments are conducted in seven steps at $t = \lbrace t_1, t_2, ..., t_7 \rbrace$ when the specimen assumes the associated configurations shown in~\cref{SetNum1} and~\cref{SetNum2}~($t_1-t_7$). Every sensing step entails in-plane harmonic excitation at a set of source points residing on the incident grid $S^{\text{inc}}$. The excitation frequency \mbox{$\omega = 72$ rad/s} is set such that the induced shear wavelength $\lambda_s$ in the specimen is approximately 0.08, giving a shear-wavelenghth-to-specimen-thickness ratio of about 4. In this setting, the phase error committed by the plane-stress approximation for the wave motion is less than $3\%$~\cite{pour2018(2)}. The incident wave interacts with both pre-existing and newly born scatterers at each $t_\alpha$ giving rise to the scattered field $\bv_\alpha$, governed by~\eqref{GEi} in Setup I -- whose pattern $\bv_\alpha{\!\!\!\obs}\!$ over the observation grid $S^{\text{obs}}$ is then computed. It must be mentioned that since the scatterers are buried within the plate thickness -- see the sectional views in~\cref{SetNum1} and~\cref{SetNum2}, our simulations are performed in three dimensions via an elastodynamic code rooted in the boundary element method~\cite{Bon1999, Fatemeh2016}. However, the nontrivial components of the computed scattered fields lay in \mbox{the $x-y$} plane, as expected in light of the earlier remarks. To study the sensitivity of evolution indictors to sensing arrangement, the incident/observation grid in \emph{Setup I} is a circle of radius $1.45$ in the mid-section of the plate, while the support of $S^{\text{inc/\nxs obs}}$ in \emph{Setup II} is the external boundary of the specimen i.e.,~a square of size $3$ $\!\times\!$ $3$. 
\vspace*{-2mm}
\subsection{Data Inversion}

With the preceding data, one may generate the evolution indicator maps affiliated with~\eqref{EIFn} in three steps, namely by:~(1) constructing the discrete scattering operators ${\text{\bf F}}_\alpha$ and ${\text{\bf F}}_\alpha^\delta$ from synthetic data related to every sensing step $t_\alpha$, (2) computing the trial signature patterns $\bPhi_L$ pertinent to a finite host domain, and (3) evaluating the differential evolution indicator in the sampling area through careful minimization of the discretized cost functional~\eqref{GCfn} as elucidated in the sequel. 

\subsubsection*{Step 1:~construction of the discrete scattering operator}
For both illumination and sensing purposes, $S^{\text{inc/\nxs obs}}$ is sampled by a uniform grid of $N$ excitation and observation points. Given the in-plane nature of wave motion, i.e., that the polarization amplitude of excitation $\bq$ and the nontrivial components of associated scattered fields $\bv_{\alpha}$ lay in \mbox{the $x-y$} plane of orthonormal bases $(\be_1,\be_2)$, the discretized scattering operator ${\text{\bf F}}_\alpha$ may be represented by a $2N\!\times 2N$ matrix with components    
\vspace{-1.5 mm}
\beq\lb{mat2}
\textrm{\bf{F}}_{\alpha}(2k\nxs+\nxs1\!:\!2k\nxs+\nxs2, \,2j\nxs+\nxs1\!:\!2j\nxs+\nxs2) ~=\, 
\left[\begin{array}{cccc}
\!W_{\alpha}^{11} \!&\! W_{\alpha}^{12}\!  \\*[1mm]
\!W_{\alpha}^{21} \!&\! W_{\alpha}^{22} \! \\
\end{array}\right] (\bx_j,\bxi_k),  \qquad j,k = 0,\ldots N-1,
\vspace{-1.5 mm}
\eeq
where~$W_{\alpha}^{\iota\upsilon}(\bx_j,\bxi_k)$ $(\iota,\upsilon\!=\!1,2)$ is the $\iota^{\textrm{th}}$ component of the displacement field measured at $\bxi_k \in S^{\text{obs}}$ due to a unit harmonic excitation applied at $\bx_j \in S^{\text{inc}}$ along the coordinate direction $\upsilon$ such that 
\vspace{-1.5 mm}
\beq\label{mat1}\nonumber
\tilde\bv_{\alpha}(\bxi_k) \,=\,
\!\left[\begin{array}{c}
\!\!\bv_{\alpha}\sip\be_1 \!\!\!\\*[0.25mm]
\!\!\bv_{\alpha}\sip\be_2\!\!\! \\
\end{array}\right] \!\!(\bxi_k) \,=\, 
\textrm{\bf{W}}_{\alpha} (\bx_j,\bxi_k) \times\!
\left[\begin{array}{c}
\!\!\bq\sip\be_1\!\!\! \\*[0.25mm]
\!\!\bq\sip\be_2 \!\!\!\\
\end{array}\right] \!\!(\bx_j). 
\vspace{-1.5 mm}
\eeq
Unless stated otherwise, we assume $N=500$. 

\emph{Noisy data.}~To account for the presence of noise in measurements, we consider the perturbed operators
\vspace{-1.5 mm}
\begin{equation}\label{DFN}
\textrm{\bf{F}}_{\alpha}^\delta \,\, \colon \!\!\!= \, (\boldsymbol{I} + \boldsymbol{N}_{\!\epsilon} ) \exs \textrm{\bf{F}}_{\alpha}, 
\vspace{-1.5 mm}
\end{equation}
where $\boldsymbol{I}$ is the $2N \times 2N$ identity matrix, and $\boldsymbol{N}_{\!\epsilon}$ is the noise matrix of commensurate dimension whose components are uniformly-distributed (complex) random variables in $[-\epsilon, \, \epsilon]^2$. In what follows, the measure of noise in data with reference to definition~\cref{Ns-op} is $\delta = \norms{\!\boldsymbol{N}_{\!\epsilon} \exs \textrm{\bf{F}}_{\alpha}\!} = 0.05$. 

\subsubsection*{Step 2:~A physics-based library of trial patterns}
This step aims to construct a suitable right hand side for the discretized far field equation in bounded and unbounded domains pertinent to the numerical experiments of this section and analytical developments of \cref{Prelim}, respectively. 

\emph{Unbounded domain in $\R^3$.}
In this case, the trial far-field pattern $\bPhi_L^\infty\in L^2(S^2)$ is given by~\cref{Phi-inf} indicating that (a) the right hand side is not only a function of the dislocation geometry $L$ but also a function of the trial opening displacement profile $\boldsymbol{a}$, and (b) computing $\bPhi_L^\infty$ generally requires an integration process at every sampling point $\bx_{\small \circ}$. In an unbounded domain, however, one may dispense with the integration process by considering a sufficiently localized (trial) FOD such as $\boldsymbol{a}(\by) = \delta (\by-\bx_{\small \circ}\!)|{\sf L}|^{-1} \bR\bn_{\small \circ}\nxs$. In this setting, without loss of generality, the dislocation support may be interpreted as an infinitesimal crack $L=\bx_{\small \circ}\!+\bR{\sf L}$ where $\bR$ is a unitary rotation matrix, and~{\sf L} represents a vanishing penny-shaped crack of unit normal~$\bn_{\small \circ}\nxs:=\lbrace 0,0,1\rbrace$. Thus, on denoting $\textrm{\bf{n}}=\bR\bn_\circ$,~\cref{Phi-inf} may be recast as    
\vspace{-1.5 mm}
\beq\lb{Phi-inf-num}
\bPhi_L^\infty(\hat\bxi) = -  \text{i} k_p \,  \hat\bxi \exs  \big[ \exs  \lambda+2\mu \exs (\textrm{\bf{n}} \cdot \hat\bxi)^2  \exs \big] \exs e^{-\text{i}k_p \hat\bxi \cdot \bx_{\small \circ}} \oplus
 -2 \text{i} \mu \exs k_s \,\hat\bxi \times \nxs (\textrm{\bf{n}} \times\hat\bxi)\exs   (\textrm{\bf{n}}\cdot\hat\bxi) \, e^{-\text{i} k_s \hat\bxi \cdot \bx_{\small \circ}}, \quad \hat\bxi \in \Omega, \, \bx_{\small \circ} \in \R^3.
 \vspace{-1.5 mm}
\eeq   

\emph{Bounded domain.} This case corresponds to the numerical experiments of this section where the background is an elastic plate ${P}$ of finite dimensions, bearing direct relevance to potential application of differential imaging to additive manufacturing and non-destructive evaluation where the target domain i.e., real-life specimen is  bounded. In this setting, it is straightforward to rigorously show that the associated patterns $\bPhi_L$ for a finite domain is governed by  
\vspace*{-4mm}
\beq\lb{PhiL}
\begin{aligned}
&\nabla \sip (\bC \colon \! \nabla \bPhi_L) \,+\, \rho \exs \omega^2\bPhi_L ~=~ \bzero \,\, &\text{in}& \,\, {P}\,\backslash L, \\*[0.25mm]
&\bn \cdot \bC \exs \colon \!  \nabla  \bPhi_L ~=~ \bzero  \quad &\,\text{on}& \,\, \partial P, \\*[0.25mm]
& \llbracket \bPhi_L \rrbracket ~=~ \boldsymbol{a} &\,\text{on}& \,\, L.
\end{aligned}     
\vspace*{-2mm}
\eeq

In what follows, the trial signatures $\textrm{\bf{v}}_{\bx_{\small \circ},\textrm{\bf{n}}}(\bxi_k)$ over the observation grid $\bxi_k \in S^{\text{obs}}$ are computed separately for every sampling point $\bx_{\small \circ} \in P$ by solving
\vspace*{-2mm}
\beq\lb{PhiL2}
\begin{aligned}
&\nabla \sip (\bC \colon \! \nabla \textrm{\bf{v}}_{\bx_{\small \circ},\textrm{\bf{n}}}) \,+\, \rho \exs \omega^2\textrm{\bf{v}}_{\bx_{\small \circ},\textrm{\bf{n}}} ~=~ \bzero \,\, &\text{in}& \,\, {P}\,\backslash L, \\*[0.25mm]
&\bn \cdot \bC \exs \colon \!  \nabla  \textrm{\bf{v}}_{\bx_{\small \circ},\textrm{\bf{n}}} ~=~ \delta (\bxi-\bx_{\small \circ}\!)|{\sf L}|^{-1} \bR\bn_{\small \circ}  \quad &\,\text{on}& \,\, \bx_{\small \circ}\!+\bR{\sf L}, \\*[0.25mm]
&\bn \cdot \bC \exs \colon \!  \nabla  \textrm{\bf{v}}_{\bx_{\small \circ},\textrm{\bf{n}}} ~=~ \bzero  \quad &\,\text{on}& \,\, \partial P,
\end{aligned}     
\vspace*{-2mm}
\eeq
within the same computational platform mentioned earlier using the boundary element method~\cite{Bon1999, Fatemeh2016}. On recalling~\eqref{mat2}, note that for every sensing point~$\bxi_k$, $\textrm{\bf{v}}_{\bx_{\small \circ},\textrm{\bf{n}}}$ has only two non-trivial components in \mbox{the $x-y$} plane, with orthonormal bases $(\be_1,\be_2)$, which are arranged into a $2N\!\times\!1$ vector as the following
\vspace{-1.5 mm}
\beq\lb{Phi-inf-Dnum}
\bPhi_{\bx_{\small \circ},\textrm{\bf{n}}}(2k+1\!:\!2k+2) ~=~ 
\! \left[\begin{array}{c}  \! \textrm{\bf{v}}_{\bx_{\small \circ},\textrm{\bf{n}}} \sip \exs\be_1 \!\! \\*[0.25mm]
\! \textrm{\bf{v}}_{\bx_{\small \circ},\textrm{\bf{n}}} \sip \exs\be_2 \!\!
\end{array} \right]\!\! (\bxi_k), \qquad k=0,\ldots N-1.
\vspace{-1.5 mm}
\eeq

\noindent \emph{Sampling.}~With reference to~\cref{SetNum1} and~\cref{SetNum2}, the search area i.e.,~the sampling region is a square $[-0.8,0.8]^2\subset P$ probed by a uniform $100 \!\times\! 100$ grid of sampling points~$\bx_{\small \circ}$ where the featured evolution indicator functionals are evaluated, while the unit circle -- spanning possible orientations for trial dislocation $L$-- is sampled by a $72$ grid of trial normal directions $\textrm{\bf{n}}=\bR\bn_\circ$. Accordingly, the evolution indicator map is constructed through minimizing~\eqref{GCfn} for a total of $M = 10000 \!\times\! 36$ trial pairs $(\bx_{\small \circ},\textrm{\bf{n}})$.     

 \begin{rem}
 It is worth mentioning that the scattering operators $\textrm{\bf{F}}_{\alpha}^\delta$ -- constructed from the forward scattering simulations of Step 1 at every sensing step $t_\alpha$, is independent of any particular choice of $L(\bx_{\small \circ},\textrm{\bf{n}})$, and thus, remain the same for all $M$ variations of $\bPhi_{\bx_{\small \circ},\textrm{\bf{n}}}$. Moreover, the right hand side of the scattering equation
 \vspace{-1.5 mm}
\beq\lb{Dff}
\textrm{\bf{F}}_{\alpha}^\delta \, \bg^{\alpha,\delta}_{\bx_{\small \circ},\textrm{\bf{n}}}~=~\bPhi_{\bx_{\small \circ},\textrm{\bf{n}}}, 
\vspace{-1.5 mm}
\eeq 
is invariant with respect to the sensing step $t_\alpha$. Therefore, for computational efficiency, one may construct a $2N \!\times\! M$ matrix that may be interpreted as a library of physically admissible patterns as the right hand side of \cref{FF} -- encompassing all choices of $L(\bx_{\small \circ},\textrm{\bf{n}})$, and solve only one equation to construct the entire imaging indicator map at every $t_\alpha$.       
 \end{rem}

\subsubsection*{Step 3:~Differential indicators of evolution}

A critical observation is that the scattering equation~\cref{Dff} is highly ill-posed at all sensing steps in that~$\,\text{det}(\textrm{\bf{F}}_{\alpha}^\delta \nxs)=0$. This problem may arise from~{(a)}~highly nonlinear nature of the inverse problem,~{(b)}~limited incident and/or ``viewing" aperture furnished by~$S^{\text{inc}\nxs /\text{obs}}$, and~{(c)}~the emergence of local (e.g., interfacial) scattered waves -- propagating in a neighborhood of certain scatterers and boundaries~\cite{Pyr1987} -- whose footprint cannot be sensed on~$S\obs$. Accordingly, \cref{Dff} will be solved via a careful regularization process via the cost functional~\eqref{GCfn}. In this setting, it is rigorously shown in \cref{Prelim} that cost functionals of type $J_{\alpha}^{\delta,\gamma}$ are convex and their minimizer can be obtained without iteration. In this vein, the discretized minimizer $\bg^{\alpha,\delta}_{\bx_{\small \circ},\textrm{\bf{n}}}$ of~\eqref{GCfn} at each sensing step $t_\alpha$ is computed by by invoking~\eqref{mat2},~\eqref{Phi-inf-Dnum}, and~\eqref{Dff} via  
\vspace{-1 mm}
\beq \lb{min-DRJ} 
\Big( \textrm{\bf{F}}_\alpha^{\delta *}\exs\textrm{\bf{F}}_\alpha^\delta  + \gamma^{\alpha,\delta}_{\bx_{\small \circ},\textrm{\bf{n}}} \exs  (\textrm{\bf{F}}_{\alpha_\sharp}^\delta)^{\nxs\frac{1}{2}*} (\textrm{\bf{F}}_{\alpha_\sharp}^\delta)^{\nxs\frac{1}{2}} + \delta \exs \gamma^{\alpha,\delta}_{\bx_{\small \circ},\textrm{\bf{n}}}  \exs \boldsymbol{I}_{2N\!\times\nxs 2N} \Big) \exs \bg^{\alpha,\delta}_{\bx_{\small \circ},\textrm{\bf{n}}}  ~=~  \textrm{\bf{F}}_\alpha^{\delta *} \bPhi_{\bx_{\small \circ},\textrm{\bf{n}}},
\vspace{-1 mm}
\eeq
where $(\cdot)^*$ is the Hermitian operator; $\textrm{\bf{F}}_{\alpha_\sharp}^\delta$ is evaluated on the basis of definitions~\cref{Fs} and~\cref{mat2}; and, following~\cite{Fatemeh2016},   
\vspace{-1 mm}  
\beq\lb{Alph}
\gamma^{\alpha,\delta}_{\bx_{\small \circ},\textrm{\bf{n}}} \,\, \colon \!\!\! = \,\, \frac{\eta^{\alpha,\delta}_{\bx_{\small \circ},\textrm{\bf{n}}}}{\norms{\textrm{\bf{F}}_\alpha^\delta\!} + \,\, \delta}.
\vspace{-1 mm}
\eeq
Here $\eta^{\alpha,\delta}_{\bx_{\small \circ},\textrm{\bf{n}}}$ is a regularization parameter computed via the Morozov discrepancy principle~\cite{Kress1999}. As a result, $\bg^{\alpha,\delta}_{\bx_{\small \circ},\textrm{\bf{n}}}$ is a $2N\!\times\! 1$ vector (or $2N\!\times\! M$ matrix for all the constructed right hand sides) identifying the structure of source densities at sensing step $t_{\alpha}$. On repeating~\eqref{min-DRJ} for all sensing steps i.e.,~$\alpha = \lbrace \circ, 1, 2, ... \rbrace$, one obtains all the arguments needed to construct a the differential evolution indicator maps.   

Next, one may compute the invariant functionals
\vspace{-1mm}
\beq \lb{Inv2N}
\begin{aligned}
\boldsymbol{\Lambda}^{\alpha,\delta}(\bg^{\alpha,\delta}_{\bx_{\small \circ},\textrm{\bf{n}}}, \bg^{\alpha+1,\delta}_{\bx_{\small \circ},\textrm{\bf{n}}}) &= \big(\, \bg^{\alpha+1,\delta}_{\bx_{\small \circ},\textrm{\bf{n}}}-\bg^{\alpha,\delta}_{\bx_{\small \circ},\textrm{\bf{n}}},\, \textrm{\bf{F}}^\delta_{\alpha_\sharp}(\bg^{\alpha+1,\delta}_{\bx_{\small \circ},\textrm{\bf{n}}}-\bg^{\alpha,\delta}_{\bx_{\small \circ},\textrm{\bf{n}}})\big) + \exs\delta \norms{\bg^{\alpha+1,\delta}_{\bx_{\small \circ},\textrm{\bf{n}}}-\bg^{\alpha,\delta}_{\bx_{\small \circ},\textrm{\bf{n}}}\!}^2, \\*[1mm]
\boldsymbol{\Upsilon}^{\alpha,\delta}(\bg^{\alpha,\delta}_{\bx_{\small \circ},\textrm{\bf{n}}}, \bg^{\alpha+1,\delta}_{\bx_{\small \circ},\textrm{\bf{n}}}) &= \big| \boldsymbol{\Lambda}^{\alpha+1,\delta}(\bzero,\bg^{\alpha+1,\delta}_{\bx_{\small \circ},\textrm{\bf{n}}}) -\boldsymbol{\Lambda}^{\alpha,\delta}(\bg^{\alpha,\delta}_{\bx_{\small \circ},\textrm{\bf{n}}}, \bzero) \big|.
\end{aligned}
\vspace{-1mm}
\eeq
Whereby, the differential evolution indicators may be computed as follows
\vspace*{-1mm}
\beq\lb{EIFnC}
\begin{aligned}
& \text{\bf I}_{\bx_{\small \circ},\textrm{\bf{n}}}^\mathcal{D^{\alpha,\delta}}(\bg^{\alpha,\delta}_{\bx_{\small \circ},\textrm{\bf{n}}}, \bg^{\alpha+1,\delta}_{\bx_{\small \circ},\textrm{\bf{n}}}) = \frac{1}{\sqrt{\boldsymbol{\Lambda}^{\alpha+1,\delta}(\bzero,\bg^{\alpha+1,\delta}_{\bx_{\small \circ},\textrm{\bf{n}}}) \big{[}1+ \boldsymbol{\Lambda}^{\alpha+1,\delta}(\bzero,\bg^{\alpha+1,\delta}_{\bx_{\small \circ},\textrm{\bf{n}}})\text{\bf{D}}_{\alpha,\delta}^{-1}(\bg^{\alpha,\delta}_{\bx_{\small \circ},\textrm{\bf{n}}}, \bg^{\alpha+1,\delta}_{\bx_{\small \circ},\textrm{\bf{n}}})\big{]}}},  \\*[1mm]
&\hat{\text{\bf I}}_{\bx_{\small \circ},\textrm{\bf{n}}}^\mathcal{D^{\alpha,\delta}}(\bg^{\alpha,\delta}_{\bx_{\small \circ},\textrm{\bf{n}}}, \bg^{\alpha+1,\delta}_{\bx_{\small \circ},\textrm{\bf{n}}}) = \frac{1}{\sqrt{\boldsymbol{\Lambda}^{\alpha,\delta}(\bg^{\alpha,\delta}_{\bx_{\small \circ},\textrm{\bf{n}}}, \bzero) + \boldsymbol{\Lambda}^{\alpha+1,\delta}(\bzero,\bg^{\alpha+1,\delta}_{\bx_{\small \circ},\textrm{\bf{n}}}) \big{[}1+ \boldsymbol{\Lambda}^{\alpha,\delta}(\bg^{\alpha,\delta}_{\bx_{\small \circ},\textrm{\bf{n}}}, \bzero)\text{\bf{D}}_{\alpha,\delta}^{{-1}}(\bg^{\alpha,\delta}_{\bx_{\small \circ},\textrm{\bf{n}}}, \bg^{\alpha+1,\delta}_{\bx_{\small \circ},\textrm{\bf{n}}})\big{]}}},
\end{aligned}
\vspace*{-1mm}
\eeq
where $\text{\bf{D}}_{\alpha,\delta} \in \lbrace \boldsymbol{\Lambda}^{\alpha,\delta}, \boldsymbol{\Upsilon}^{\alpha,\delta} \rbrace$. Upon introducing
\beq
(\bg^{\alpha,\delta}_{\bx_{\small \circ}}, \bg^{\alpha+1,\delta}_{\bx_{\small \circ}}) \,\,\colon \!\!= \,\, \text{argmin}_{(\bg^{\alpha,\delta}_{\bx_{\small \circ},\textrm{\bf{n}}}, \bg^{\alpha+1,\delta}_{\bx_{\small \circ},\textrm{\bf{n}}})} \big( \text{\bf I}_{\bx_{\small \circ},\textrm{\bf{n}}}^\mathcal{D^{\alpha,\delta}} \big), \quad (\hat{\bg}^{\alpha,\delta}_{\bx_{\small \circ}}, \hat{\bg}^{\alpha+1,\delta}_{\bx_{\small \circ}}) \,\,\colon \!\!= \,\, \text{argmin}_{(\bg^{\alpha,\delta}_{\bx_{\small \circ},\textrm{\bf{n}}}, \bg^{\alpha+1,\delta}_{\bx_{\small \circ},\textrm{\bf{n}}})} \big( \hat{\text{\bf I}}_{\bx_{\small \circ},\textrm{\bf{n}}}^\mathcal{D^{\alpha,\delta}} \big),
\eeq
one obtains the evolution indicator maps
\vspace*{-1.5mm}
\beq\lb{EIFN}
\begin{aligned}
& \text{\bf I}_{\bx_{\small \circ}}^\mathcal{D^{\alpha,\delta}}(\bg^{\alpha,\delta}_{\bx_{\small \circ}}, \bg^{\alpha+1,\delta}_{\bx_{\small \circ}}) = \frac{1}{\sqrt{\boldsymbol{\Lambda}^{\alpha+1,\delta}(\bzero,\bg^{\alpha+1,\delta}_{\bx_{\small \circ}}) \big{[}1+ \boldsymbol{\Lambda}^{\alpha+1,\delta}(\bzero,\bg^{\alpha+1,\delta}_{\bx_{\small \circ}})\text{\bf{D}}_{\alpha,\delta}^{-1}(\bg^{\alpha,\delta}_{\bx_{\small \circ}}, \bg^{\alpha+1,\delta}_{\bx_{\small \circ}})\big{]}}},  \\*[0.25mm]
&\hat{\text{\bf I}}_{\bx_{\small \circ}}^\mathcal{D^{\alpha,\delta}}(\bg^{\alpha,\delta}_{\bx_{\small \circ}}, \bg^{\alpha+1,\delta}_{\bx_{\small \circ}}) = \frac{1}{\sqrt{\boldsymbol{\Lambda}^{\alpha,\delta}(\bg^{\alpha,\delta}_{\bx_{\small \circ}}, \bzero) + \boldsymbol{\Lambda}^{\alpha+1,\delta}(\bzero,\bg^{\alpha+1,\delta}_{\bx_{\small \circ}}) \big{[}1+ \boldsymbol{\Lambda}^{\alpha,\delta}(\bg^{\alpha,\delta}_{\bx_{\small \circ}}, \bzero)\text{\bf{D}}_{\alpha,\delta}^{{-1}}(\bg^{\alpha,\delta}_{\bx_{\small \circ}}, \bg^{\alpha+1,\delta}_{\bx_{\small \circ}})\big{]}}},
\end{aligned}
\vspace*{-1.5mm}
\eeq

Here, $\text{\bf I}_{\bx_{\small \circ}}^\mathcal{D^{\alpha,\delta}}\!\!$ and $\hat{\text{\bf I}}_{\bx_{\small \circ}}^\mathcal{D^{\alpha,\delta}}\!\!$ canvas the support of  evolution between two sensing steps $\alpha$ and $\alpha+1$. More specifically, $\hat{\text{\bf I}}_{\bx_{\small \circ}}^\mathcal{D^{\alpha,\delta}}\!\!$ assumes~{(a)}~its highest values at those sampling points that meet the newly developed or evolved scatterers between $t_\alpha$ and $t_{\alpha+1}$, and~{\bf (b)}~near zero values everywhere else including the loci of (unknown) stationary scatterers within $[t_\alpha\,\, t_{\alpha+1}]$ timeframe.   

\begin{rem}
To gain better insight into the effectiveness of the proposed approach, the reconstructed evolution indicators will be compared to their corresponding~GLSM map. With reference to~\cref{GLSMgs}, the latter can be computed via 
\vspace*{-1.5mm}
\beq\lb{GLSM-Dgs}
\text{\bf I}^{\mathcal{G}_\sharp^{\alpha,\delta}}_{\bx_{\small \circ}} \,\, = \,\, \dfrac{1}{\sqrt{\norms{\!(\textrm{\bf{F}}^\delta_\sharp)^{\frac{1}{2}} \exs \bg^{\alpha,\delta}_{\bx_{\small \circ}} \nxs}^2 \exs+\,\, \delta \! \norms{\bg^{\alpha,\delta}_{\bx_{\small \circ}} \!}^2}}, \qquad 
\textcolor{black}{
\bg^{\alpha,\delta}_{\bx_{\small \circ}} \,\,\colon \!\!= \,\, \text{argmin}_{\bg^{\alpha,\delta}_{\bx_{\small \circ},\textrm{\bf{n}}}} \norms{\bg^{\alpha,\delta}_{\bx_{\small \circ},\textrm{\bf{n}}}}^2_{L^2}.} 
\vspace*{-1.5mm}
\eeq  
\end{rem}

\subsection{Simulation results}\label{CompR}

The synthetic scattered data deployed to generated the ensuing results are perturbed by $\delta \!=\! 5\%$ white noise.~\cref{RD1} illustrates the full-aperture GLSM reconstruction of a progressive damage zone using scattered field data computed in seven sensing steps $t _{\alpha}\!=\! \lbrace t_1, t_2, ..., t_7 \rbrace$. Every panel in~\cref{RD1} demonstrates the distribution of indicator functional $\text{\bf I}^{\mathcal{G}_\sharp^{\alpha,\delta}}_{\bx_{\small \circ}}\!\!$ over the indicated sampling area. Note that the GLSM functional at every $t_\alpha$ is solely dependent upon the far field data at that step $\bv_{\alpha}{\!\!\!\obs}\!$ i.e.,~this imaging indicator does not require knowledge of the sequential sensory data ($\bv_{\alpha+1}{\!\!\!\!\!\!\!\!\!\!\obs}\,\,$ or $\bv_{\alpha-1}{\!\!\!\!\!\!\!\!\!\!\obs}\,\,$).~\cref{RA1} shows parallel results pertaining to an evolving bubble zone related to the second set numerical experiments depicted in~\cref{SetNum2}.   

\cref{ED1} illustrates the sequential reconstruction of newborn fractures emerging between every pair of successive sensing steps $t_{\alpha}\!\nxs-\! t_{\alpha+1}$ for $\alpha = \lbrace 1, 2, ..., 6 \rbrace$ by way of the proposed differential evolution indicators. Each panel in~\cref{ED1} demonstrates the distribution of evolution indicator functional ${\text{\bf I}_{\bx_{\small \circ}}^\mathcal{D^{\alpha,\delta}}\!\!}$ over the sampling area. According to \eqref{EIFN}, evaluation of the differential indicators requires pairs of sensory data in the form of $(\bv_{\alpha}{\!\!\!\obs}\!,\bv_{\alpha+1}{\!\!\!\!\!\!\!\!\!\!\obs}\,\,)$.~\cref{EA1} shows similar plots corresponding to an evolving bubble zone.
 
\cref{CD} and~\cref{CA} compare the performance of existing GLSM criteria with that of the proposed differential evolution indicators in reconstructing highly scattering damage and bubble zones where the distance between scatterers is of the order of a fraction of the shear wavelength. Here, the GLSM maps of damage/bubble zone at $t_7$ are depicted against that assembled map of fractures (\emph{resp.} pores) constructed by averaging the differential indictor maps shown in~\cref{ED1} (\emph{resp.}~\cref{EA1}). Note that the pre-existing bubbles shown in~\cref{SetNum2}~($t_1$), and reconstructed in~\cref{RA1}~($t_1$), do not appear in the evolution indicator map of~\cref{CA}~(b) depicting the reconstructed evolution of pore zone between $t_1\!\nxs-\!t_7$. 

\section{Conclusion}\lb{Conc}   

A robust framework for waveform tomography of {\emph{progressive evolution}} in highly scattering materials of uncertain structure is developed. Our three-tier imaging platform is inherently non-iterative enabling fast inversion of dense data in support of the real-time sensing. In addition, this method allows for concurrent reconstruction of multiple interfacial and volumetric scatterers of arbitrary geometry. This imaging solution enables targeted characterization of active zones with little sensitivity to the noise in measurements while remaining systematically agnostic with respect to uncertainties of the host domain. As a perspective it would be interesting to incorporate a multi-frequency framework that would enable multiscale characterization of evolution without the need to reconstruct the entire domain across pertinent scales which may be practically insurmountable. The proposed imaging solution is formulated in a generic framework whose validity is rigorously established. As a result, this method may potentially contribute to a wide range of civil, aerospace and mechanical systems.

\begin{figure}[!h]
\center\includegraphics[width=0.8\linewidth]{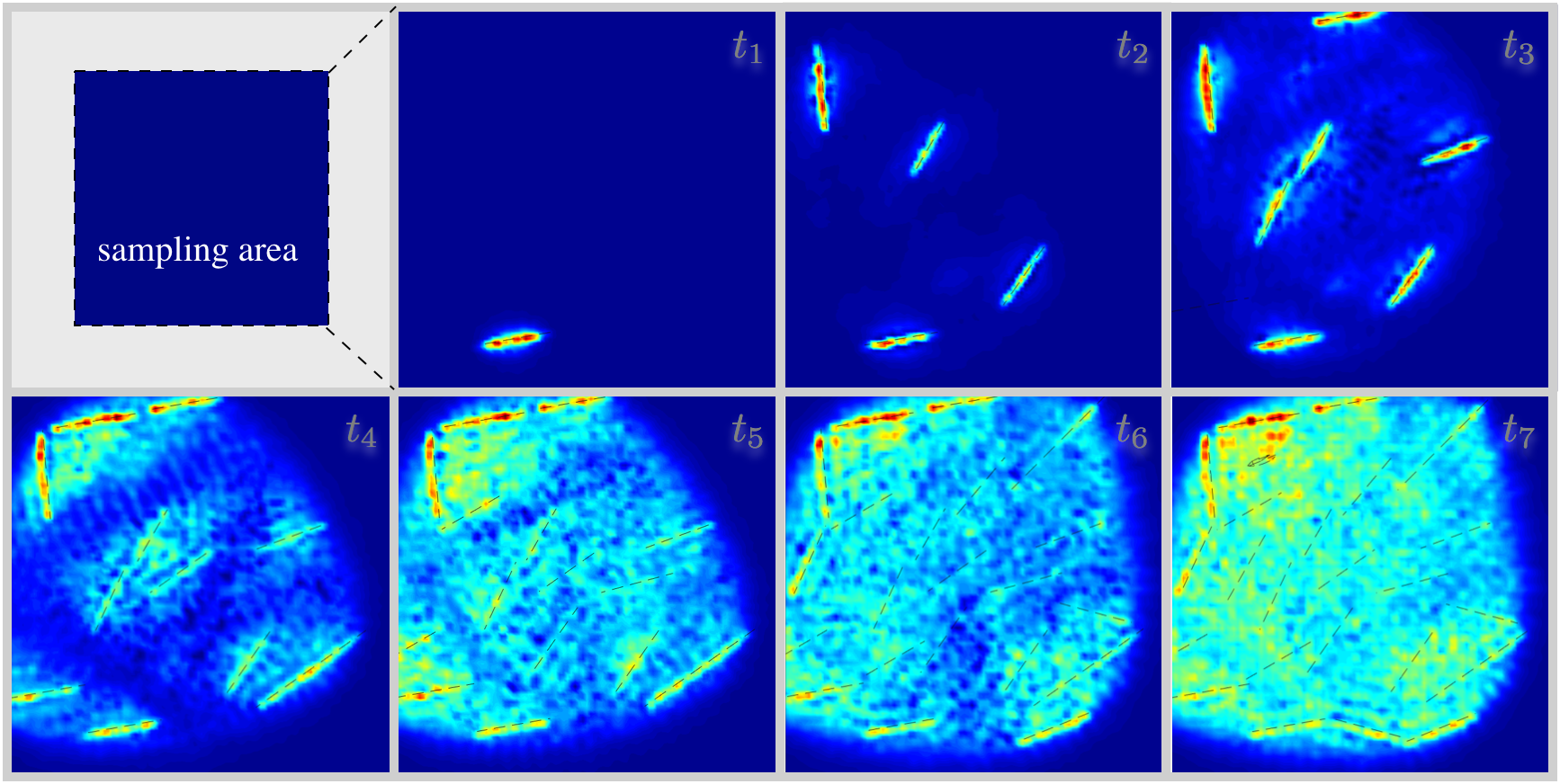} \vspace*{-3mm} 
\caption{Disjoint reconstruction of an evolving damage zone via the existing indicator {\small ${\text{\bf I}^{\mathcal{G}_\sharp^{\alpha,\delta}}_{\bx_{\small \circ}}}\!\!$} at time steps $t _{\alpha}\!=\! \lbrace t_1, t_2, ..., t_7 \rbrace$. Each panel is a GLSM map constructed via \eqref{GLSM-Dgs} on the basis of scattered field data $\bv_{\alpha}{\!\!\!\obs}\!$ exclusively captured at the indicated sensing step $t_\alpha$.} \lb{RD1}
\vspace*{-2mm}
\end{figure} 

\begin{figure}[!h]
\center\includegraphics[width=0.62\linewidth]{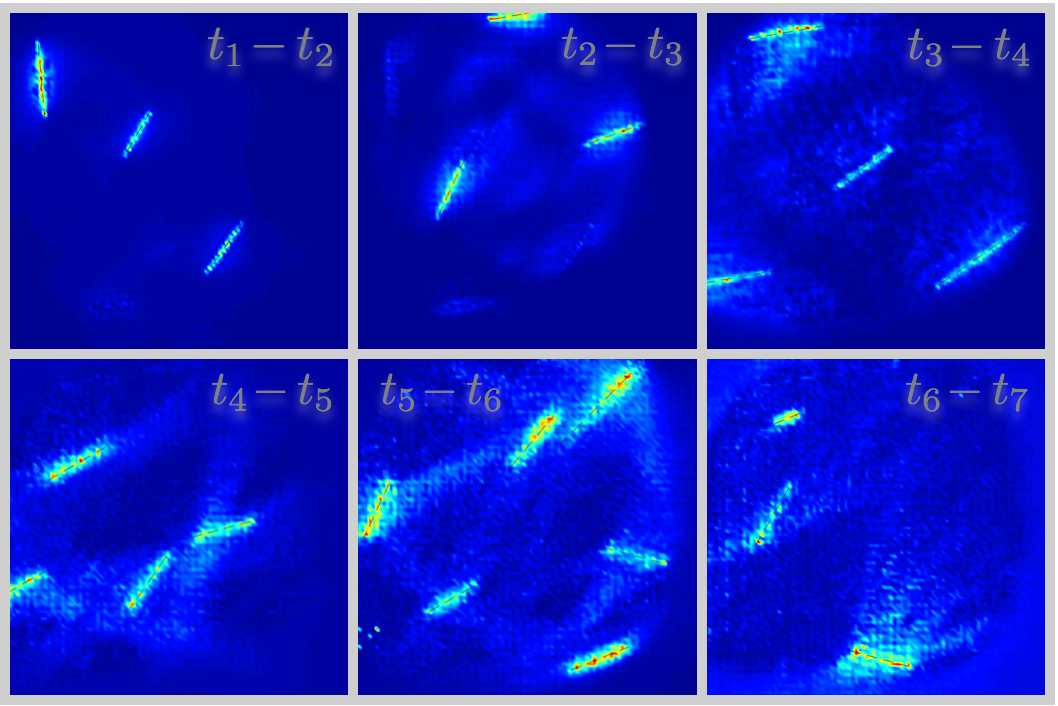} \vspace*{-3mm} 
\caption{Differential reconstruction of newborn fractures, via the proposed indicator ${\text{\bf I}_{\bx_{\small \circ}}^\mathcal{D^{\alpha,\delta}}\!\!}$, within every sensing sequence $t_{\alpha}\!\nxs-\! t_{\alpha+1}$ for $\alpha = \lbrace 1, 2, ..., 6 \rbrace$. Each panel is a Differential Evolution Indicator map constructed via \eqref{EIFN} on the basis of sequential scattered field data $(\bv_{\alpha}{\!\!\!\obs}\!,\bv_{\alpha+1}{\!\!\!\!\!\!\!\!\!\!\obs}\,\,)$, captured at the indicated sensing steps $t_\alpha, t_{\alpha+1}$.} \lb{ED1}
\vspace*{-4mm}
\end{figure} 

\begin{figure}[!h]
\center\includegraphics[width=0.62\linewidth]{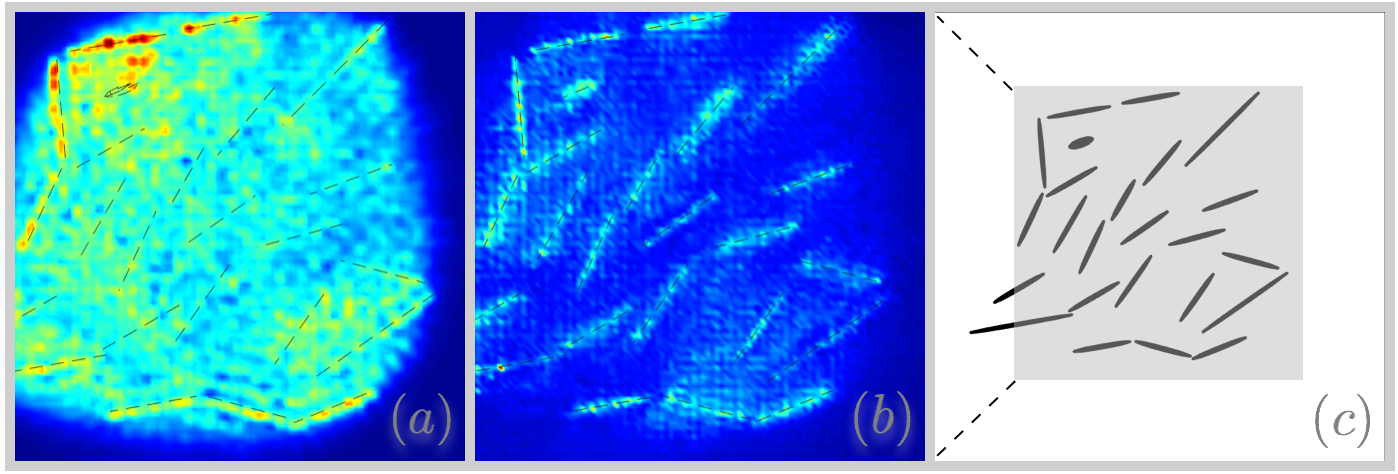} \vspace*{-3mm} 
\caption{GLSM~vs.~Differential Evolution Indicators:~(a)~GLSM map of damage zone within the sampling area at $t_7$, (b) assembled map of fractures constructed by averaging the differential indictor maps shown in~\cref{EA1}, and (c) true anatomy of the damage zone where the shaded area highlights the sampling region.} \lb{CD}
\vspace*{-5mm}
\end{figure} 

\pagebreak

\begin{figure}[!h]
\center\includegraphics[width=0.8\linewidth]{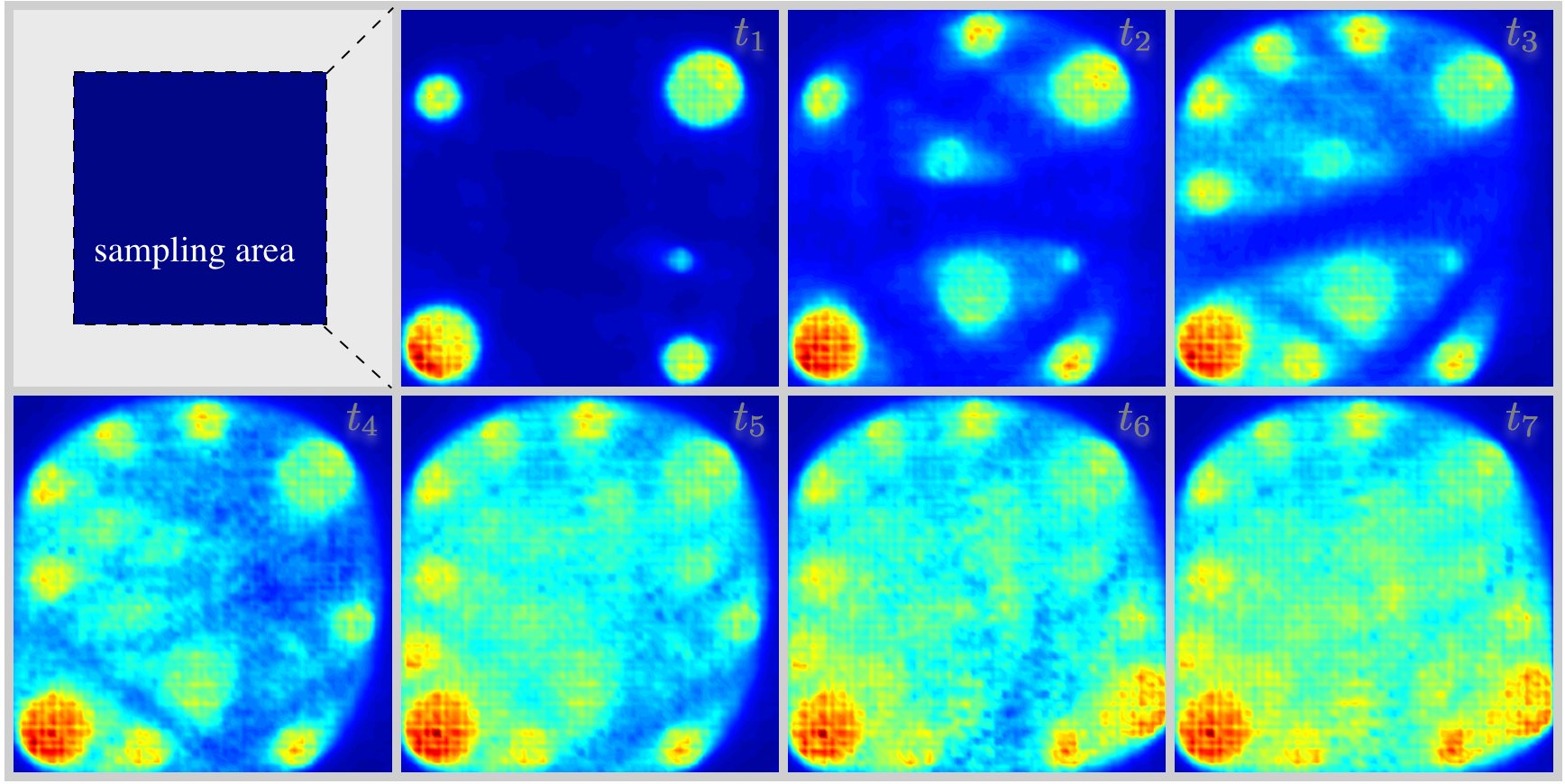} \vspace*{-3mm} 
\caption{Disjoint reconstruction of an evolving bubble zone via the existing indicator {\small ${\text{\bf I}^{\mathcal{G}_\sharp^{\alpha,\delta}}_{\bx_{\small \circ}}}\!\!$} at time steps $t _{\alpha}\!=\! \lbrace t_1, t_2, ..., t_7 \rbrace$. Each panel is a GLSM map constructed via \eqref{GLSM-Dgs} on the basis of scattered field data $\bv_{\alpha}{\!\!\!\obs}\!$ exclusively captured at the indicated sensing step $t_\alpha$.} \lb{RA1}
\vspace*{-5mm}
\end{figure} 

\begin{figure}[!h]
\center\includegraphics[width=0.62\linewidth]{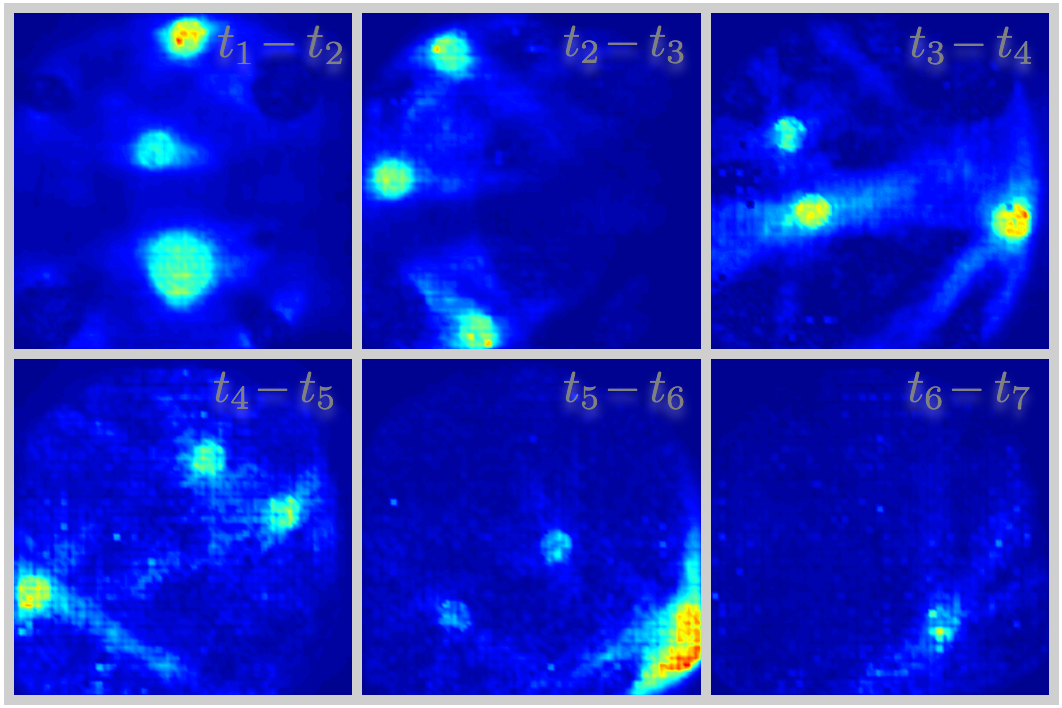} \vspace*{-3mm} 
\caption{Differential reconstruction of newborn pores, via the proposed indicator ${\text{\bf I}_{\bx_{\small \circ}}^\mathcal{D^{\alpha,\delta}}\!\!}$, within every sensing sequence $t_{\alpha}\!\nxs-\! t_{\alpha+1}$ for $\alpha = \lbrace 1, 2, ..., 6 \rbrace$. Each panel is a Differential Evolution Indicator map constructed via \eqref{EIFN} on the basis of sequential scattered field data $(\bv_{\alpha}{\!\!\!\obs}\!,\bv_{\alpha+1}{\!\!\!\!\!\!\!\!\!\!\obs}\,\,)$, captured at the indicated sensing steps $t_\alpha, t_{\alpha+1}$.} \lb{EA1}
\vspace*{-4mm}
\end{figure} 

\begin{figure}[!h]
\center\includegraphics[width=0.62\linewidth]{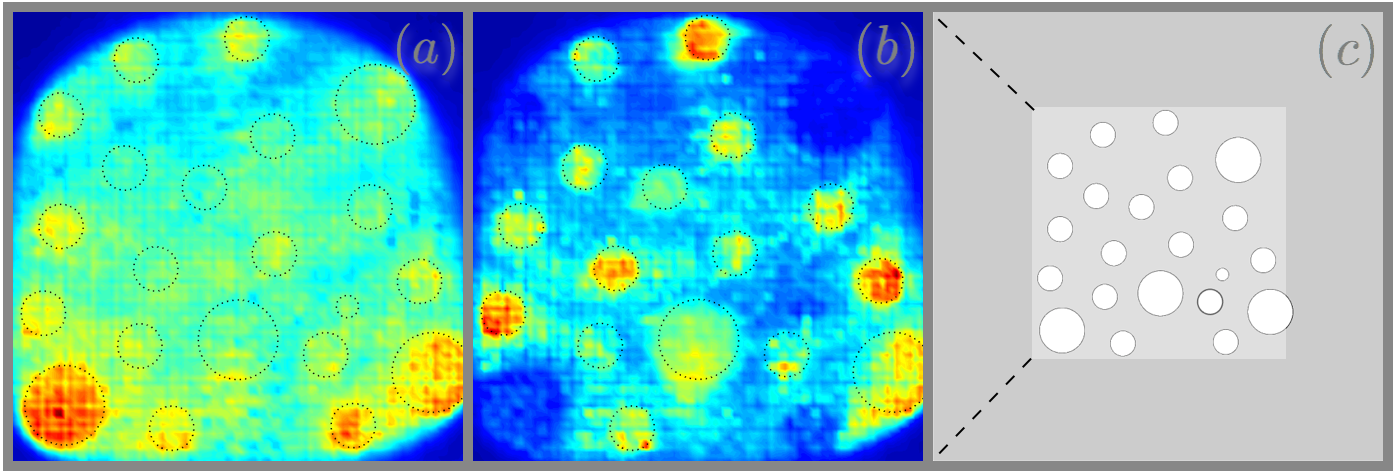} \vspace*{-3mm} 
\caption{GLSM~vs.~Differential Evolution Indicators:~(a)~GLSM map of bubble zone within the sampling area at $t_7$, (b) assembled map of pores constructed by averaging the differential indictor maps shown in~\cref{EA1}, and (c) true anatomy of the pore zone where the shaded area highlights the sampling region.} \lb{CA}
\vspace*{-5mm}
\end{figure}

\section{Acknowledgements}\lb{Ac}   

The corresponding author kindly acknowledges the support provided by the University of Colorado Boulder through the Imaging Science seed grants and FP's startup. This work utilized the RMACC Summit supercomputer, which is supported by the National Science Foundation (awards ACI-1532235 and ACI-1532236), the University of Colorado Boulder, and Colorado State University. The Summit supercomputer is a joint effort of the University of Colorado Boulder and Colorado State University.

\bibliographystyle{siamplain}
\bibliography{inverse}

\end{document}